# Lattice point problems and distribution of values of quadratic forms

By V. Bentkus and F. Götze*

## Abstract

For $d$-dimensional *irrational* ellipsoids $E$ with $d \geq 9$ we show that the number of lattice points in $rE$ is approximated by the volume of $rE$, as $r$ tends to infinity, up to an error of order $o(r^{d-2})$. The estimate refines an earlier authors' bound of order $\mathcal{O}(r^{d-2})$ which holds for *arbitrary* ellipsoids, and is optimal for rational ellipsoids. As an application we prove a conjecture of Davenport and Lewis that the gaps between successive values, say $s < n(s)$, $s, n(s) \in Q[\mathbb{Z}^d]$, of a positive definite irrational quadratic form $Q[x]$, $x \in \mathbb{R}^d$, are shrinking, i.e., that $n(s) - s \to 0$ as $s \to \infty$, for $d \geq 9$. For comparison note that $\sup_s(n(s) - s) < \infty$ and $\inf_s(n(s) - s) > 0$, for rational $Q[x]$ and $d \geq 5$. As a corollary we derive Oppenheim's conjecture for indefinite irrational quadratic forms, i.e., the set $Q[\mathbb{Z}^d]$ is dense in $\mathbb{R}$, for $d \geq 9$, which was proved for $d \geq 3$ by G. Margulis [Mar1] in 1986 using other methods. Finally, we provide explicit bounds for errors in terms of certain characteristics of trigonometric sums.

## 1. Introduction and results

Let $\mathbb{R}^d$, $1 \leq d < \infty$, denote a real $d$-dimensional Euclidean space with scalar product $\langle \cdot, \cdot \rangle$ and norm

$$|x|^2 = \langle x, x \rangle = x_1^2 + \cdots + x_d^2, \qquad \text{for } x = (x_1, \ldots, x_d) \in \mathbb{R}^d.$$

We shall use as well the norms $|x|_1 = \sum_{j=1}^d |x_j|$ and $|x|_\infty = \max\{ |x_j| : 1 \leq j \leq d\}$. Let $\mathbb{Z}^d$ be the standard lattice of points with integer coordinates in $\mathbb{R}^d$.

---

*Research supported by the SFB 343 in Bielefeld.

1991 *Mathematics Subject Classification.* 11P21.

*Key words and phrases.* lattice points, ellipsoids, rational and irrational quadratic forms, positive and indefinite quadratic forms, distribution of values of quadratic forms, Oppenheim conjecture, Davenport-Lewis conjecture.

For a (measurable) set $B \subset \mathbb{R}^d$, let vol $B$ denote the Lebesgue measure of $B$, and let $\mathrm{vol}_{\mathbb{Z}}\, B$ denote the lattice volume of $B$, that is the number of points in $B \cap \mathbb{Z}^d$.

Consider a quadratic form,

$$Q[x] \overset{\mathrm{def}}{=} \langle Qx, x\rangle, \qquad \text{for} \quad x \in \mathbb{R}^d,$$

where $Q : \mathbb{R}^d \to \mathbb{R}^d$ denotes a symmetric linear operator with nonzero eigenvalues, say $q_1, \dots, q_d$. Write

$$q_0 = \min_{1\le j\le d} |q_j|, \qquad q = \max_{1\le j\le d} |q_j|. \tag{1.1}$$

We assume that the form in nondegenerate, that is, that $q_0 > 0$. Thus, without loss of generality we can and shall assume throughout that $q_0 = 1$, and hence $q \ge 1$.

Define the sets

$$E_s = \{x \in \mathbb{R}^d : \ Q[x] \le s\}, \qquad \text{for } s \in \mathbb{R}.$$

If the operator $Q$ is positive definite (henceforth called briefly positive), that is, $Q[x] > 0$, for $x \ne 0$, then $E_s$ is an ellipsoid.

Recall that a quadratic form $Q[x]$ with a nonzero matrix $Q = (q_{ij})$, $1 \le i$, $j \le d$, is rational if there exists an $M \in \mathbb{R}$, $M \ne 0$, such that the matrix $MQ$ has integer entries only; otherwise it is called irrational. We identify the matrix of $Q[x]$ with the operator $Q$.

Our main result, Theorems 2.1 and 2.2, yields the following Theorems 1.1, 1.3 and 1.5 and Corollaries 1.2, 1.4, 1.6 and 1.7, proofs of which we provide in Section 2.

THEOREM 1.1. *Assume that $Q$ is positive and $d \ge 9$. Then*

$$\sup_{a\in\mathbb{R}^d} \Delta(s, Q, a) \overset{\mathrm{def}}{=} \sup_{a\in\mathbb{R}^d} \left| \frac{\mathrm{vol}_{\mathbb{Z}}\,(E_s + a) - \mathrm{vol}\, E_s}{\mathrm{vol}\, E_s} \right| = o(s^{-1}), \tag{1.2}$$

*as $s \to \infty$ if and only if $Q$ is irrational.*

The estimate of Theorem 1.1 refines an explicit bound of order $\mathcal{O}(s^{-1})$ obtained by the authors (henceforth called [BG1]) for *arbitrary* ellipsoids. That result has been proved using probabilistic notions and a version of the basic inequality (see (3.12) below) for trigonometric sums. Some methods of that proof will be used again in this paper. An alternative proof using more extensively the method of large sieves appeared as [BG4]. In the case of *rational* ellipsoids the bound $\mathcal{O}(s^{-1})$ is optimal. For arbitrary ellipsoids Landau [La1] obtained the estimate $\mathcal{O}(s^{-1+1/(1+d)})$, $d \ge 1$. This result has been extended by Hlawka [H] to convex bodies with smooth boundary and strictly positive Gaussian curvature. Hlawka's estimate has been improved by Krätzel and Nowak ([KN1], [KN2]) to $\mathcal{O}(s^{-1+\lambda})$, where $\lambda = 5/(6d+2)$, for $d \geq 8$, and $\lambda = 12/(14d+8)$, for $3 \leq d \leq 7$. For special ellipsoids a number of particular results is available. For example, the error bound $\mathcal{O}(s^{-1})$ holds for $d \geq 5$ and rational $Q$ (see Walfisz [Wa1], $d \geq 9$, and Landau [La2], $d \geq 5$). Jarnik [J1] proved the same bound for diagonal $Q$ with arbitrary (nonzero) real entries. For a discussion see the monograph Walfisz [Wa2].

Theorem 1.1 is applicable to irrational ellipsoids with arbitrary center for $d \geq 9$. It extends the bound of order $o(s^{-1})$ of Jarnik and Walfisz [JW] for diagonal irrational $Q$ of dimension $d \geq 5$. They showed that $o(s^{-1})$ is optimal, that is, for any function $\xi$ such that $\xi(s) \to \infty$, as $s \to \infty$, there exists an irrational diagonal form $Q[x]$ such that

$$\limsup_{s\to\infty}\ s\,\xi(s)\,\Delta(s,Q,0) = \infty.$$

See Theorem 1.3 for an estimate of the remainder term in (1.2) in terms of certain characteristics of trigonometric sums.

*Gaps between values of positive quadratic forms.* Let $s, n(s) \in Q[\mathbb{Z}^d]$, $s < n(s)$, denote successive values of $Q[x]$. Davenport and Lewis [DL] conjectured that the distance between successive values of the quadratic form $Q[x]$ converges to zero as $s \to \infty$, provided that the dimension $d \geq 5$ and $Q$ is irrational. Theorem 1.1 combined with Theorem 1.1 of [BG4] provides a complete solution of this problem for $d \geq 9$. Introduce the maximal gap $d(\tau; Q, a) = \sup\{n(s) - s : s \geq \tau\}$ between values $Q[x-a]$ in the interval $[\tau, \infty)$.

COROLLARY 1.2. *Assume that $d \geq 9$ and $Q[x]$ is positive definite. If the quadratic form is irrational then $\sup_{a\in\mathbb{R}^d} d(\tau; Q, a) \to 0$, as $\tau \to \infty$. If $Q$ is rational then $\sup_{\tau\geq 0} d(\tau; Q, a) < \infty$. If both $Q$ and $a$ are rational then $\inf_s \big(n(s) - s\big) > 0$.*

Answering a question of T. Esterman whether gaps must tend to zero for large dimensional positive forms, Davenport and Lewis [DL] proved the following: Assume that $d \geq d_0$ with some sufficiently large $d_0$. Let $\varepsilon > 0$. Suppose that $y \in \mathbb{Z}^d$ has a sufficiently large norm $|y|_\infty$. Then there exists $x \in \mathbb{Z}^d$ such that

$$\big|Q[y+x] - Q[x]\big| < \varepsilon. \tag{1.3}$$

Of course (1.3) does not rule out the possibility of arbitrarily large gaps between possible clusters of values $Q[x]$, $x \in \mathbb{Z}^d$. The result of [DL] was improved by Cook and Raghavan [CR]. They obtained the estimate $d_0 \leq 995$ and provided a lower bound for the number of solutions $x \in \mathbb{Z}^d$ of the inequality (1.3). See the reviews of Lewis [Le] and Margulis [Mar2].

In order to provide bounds concerning gaps between values of positive quadratic forms and lattice point approximations for ellipsoids we need additional notation. Introduce the trigonometric sum

$$\varphi_a(t;s) = \bigl(2[\sqrt{s}\,]+1\bigr)^{-3d} \Bigl|\sum \exp\bigl\{\,itQ[x_1+x_2+x_3-a]\,\bigr\}\Bigr|, \tag{1.4}$$

where the sum is taken over all $x_1, x_2, x_3 \in \mathbb{Z}^d$ such that $|x_j|_\infty \le \sqrt{s}$, for $j = 1, 2, 3$. Notice that the trigonometric sum (1.4) is normalized so that $\varphi_a(t;s) \le \varphi_a(0;s) = 1$.

Theorem 6.1 shows that for irrational $Q$ and any fixed $0 < \delta_0 \le T < \infty$ the trigonometric sum $\varphi_a$ satisfies

$$\lim_{s\to\infty} \sup_{a\in\mathbb{R}^d} \sup_{\delta_0\le t\le T} \varphi_a(t;s) = 0. \tag{1.5}$$

Simple selection arguments show that (1.5) yields that for irrational $Q$ there exist sequences $T(s) \uparrow \infty$, $T(s) \ge 1$, and $\delta_0(s) \downarrow 0$ such that

$$\lim_{s\to\infty} \sup_{a\in\mathbb{R}^d} \sup_{\delta_0(s)\le t\le T(s)} \varphi_a(t;s) = 0. \tag{1.6}$$

The relation (6.5) shows that

$$\lim_{s\to\infty} \sup_{a\in\mathbb{R}^d} \sup_{s^{-1/2}\le t\le \delta_0(s)} \varphi_a(t;s) = 0,$$

for any nondegenerate $Q$. Hence, the irrationality of $Q$ is equivalent to the following condition: there exist $T(s) \uparrow \infty$ such that

$$\lim_{s\to\infty} \gamma\bigl(s, T(s)\bigr) = 0, \quad \text{where } \gamma\bigl(s,T\bigr) = \sup_{a\in\mathbb{R}^d} \sup_{s^{-1/2}\le t\le T} \varphi_a(t;s). \tag{1.7}$$

Finally, given $d \ge 9$, $0 < \varepsilon < 1 - 8/d$ and $Q$, introduce the quantity

$$\rho(s) = s^{1-\zeta} + \frac{1}{T(s)} + \Bigl(\gamma\bigl(s, T(s)\bigr)\Bigr)^{1-8/d-\varepsilon} T^{\varepsilon}(s), \quad \zeta \overset{\text{def}}{=} \frac{1}{2}\left[\frac{d-1}{2}\right], \tag{1.8}$$

on which our estimates will depend. Without loss of generality we can assume that $T(s)$ in (1.7) and (1.8) are chosen so that

$$\lim_{s\to\infty} \rho(s) = 0, \tag{1.9}$$

for irrational $Q$. Indeed, if (1.9) does not hold, we can replace $T(s)$ in (1.7) by

$$\min\Bigl\{T(s);\ \Bigl(\gamma\bigl(s, T(s)\bigr)\Bigr)^{-(1-8/d-\varepsilon)/(2\varepsilon)}\Bigr\}.$$

We shall write $A \ll_d B$ if there exist a constant $c_d$ depending on $d$ only and such that $A \le c_d\, B$.

Theorem 1.3. *Assume that the operator $Q$ is positive, $d \geq 9$ and $0 < \varepsilon < 1 - 8/d$. Then we have*

$$\big|\operatorname{vol}_{\mathbb{Z}}(E_s + a) - \operatorname{vol} E_s\big| \ll_{d,\varepsilon} (s+1)^{d/2}\, q^d\, \rho(s)\, s^{-1}. \tag{1.10}$$

Theorem 1.1 is an immediate corollary of Theorem 1.3.

If we choose $T(s) = 1$ in (1.8) and use $\gamma\big(s, T(s)\big) \leq 1$, then (1.10) yields

$$\big|\operatorname{vol}_{\mathbb{Z}}(E_s + a) - \operatorname{vol} E_s\big| \ll_d (s+1)^{d/2}\, q^d\, s^{-1}.$$

This slightly improves the bound $(s+1)^{d/2}\, q^{d+2}\, s^{-1}$ given in [BG4].

An inspection of proofs shows that Theorem 1.3 holds for any $\zeta < d/4$. Moreover, the main result, Theorems 2.1 and 2.2, can be proved for any real $2 \leq p < d/2$ (with the expansion in (2.8) defined by the same formula as in the case $p \in \mathbb{N}$). The assumption $p \in \mathbb{N}$ is made for technical convenience only. Hence, the presence in (1.8) of the term $s^{1-\zeta}$ shows that our bound (1.10) can not decrease faster than $\mathcal{O}(s^{d/4+\delta})$, $\delta > 0$.

Write

$$\rho_0(s) = \sup_{\tau \geq s} \rho(\tau).$$

Corollary 1.4. *Assume that the operator $Q$ is positive and $d \geq 9$. Then*

$$\frac{\operatorname{vol}_{\mathbb{Z}}\Big((E_{\tau+\delta} + a) \setminus (E_\tau + a)\Big)}{\operatorname{vol}(E_{\tau+\delta} \setminus E_\tau)} = 1 + \mathcal{R}, \quad \text{for } \tau \geq s \geq 1,\ \delta > 0, \tag{1.11}$$

*where $\mathcal{R}$ satisfies $|\mathcal{R}| \ll_{d,\varepsilon} q^{3d/2}\rho_0(s)/\delta$. In particular, the maximal gap $d(s, Q, a)$ satisfies*

$$d(s, Q, a) \ll_{d,\varepsilon} q^{3d/2}\, \rho_0(s), \qquad \text{for } s \geq 1. \tag{1.12}$$

The relation (1.11) gives an estimate of the number of values of a positive quadratic form in an interval $(s, s+\delta]$, counting these values according to their multiplicities. More precisely, a value, say $\tau = Q[x-a]$, is counted $\operatorname{vol}_{\mathbb{Z}}\big\{z : \tau = Q[z-a]\big\}$ times. In the case of irrational $Q$ the approximation (1.11) may be applied for intervals of shrinking size $\delta = \delta(s) \to 0$ as $s \to \infty$. The approximation error in (1.11) still satisfies $\mathcal{R} \to 0$ for shrinking intervals such that $\delta/\rho_0(s) \to \infty$ as $s \to \infty$. The estimate (1.12) provides an upper bound for the maximal gap $d(s, Q, a)$ between values to the right of a value $s \geq 1$, for positive $Q$. In particular, we get, for irrational $Q$, $d(s, Q, a) \to 0$ as $s \to \infty$, uniformly with respect to $a$.

*The Oppenheim conjecture.* Write

$$m(Q) = \inf\Big\{\big|Q[x]\big| : x \neq 0,\ x \in \mathbb{Z}^d\Big\}.$$

Oppenheim ([O1], [O2]) conjectured that $m(Q) = 0$, for $d \geq 5$ and irrational *indefinite* $Q$. This conjecture has been extensively studied, see the review of Margulis [Mar2]. A stronger version was finally proved by Margulis [Mar1]: $m(Q) = 0$, for $d \geq 3$ and irrational indefinite $Q$. In 1953, A. Oppenheim proved in three papers that such a result is equivalent to the following: for irrational $Q$ and $d \geq 3$, the set $Q[\mathbb{Z}^d]$ is dense in $\mathbb{R}^d$. See the discussion in [Mar2, p. 284]. In particular, $d(\tau, Q, 0) \equiv 0$, for all $\tau$, which is impossible for positive forms.

The quantitative version of Oppenheim's conjecture was developed by Dani and Margulis [DM] and Eskin, Margulis and Mozes [EMM]. Let $M : \mathbb{R}^d \to [0, \infty)$ be any *continuous* function such that $M(tx) = |t| M(x)$, for all $t \in \mathbb{R}$ and $x \in \mathbb{R}^d$, and such that $M(x) = 0$ if and only if $x = 0$. The function $M$ is the Minkowski functional of the set

$$\Omega = \big\{ x \in \mathbb{R}^d : \ M(x) \leq 1 \big\}. \tag{1.13}$$

In particular, the set $\Omega$ is a star-shaped closed bounded set with the nonempty interior containing zero. For an interval $I = (\alpha, \beta]$ define the set $W = \big\{ x \in \mathbb{R}^d : Q[x - a] \in I \big\}$. Assuming that $d \geq 5$ and that the quadratic form $Q[x]$ is irrational and indefinite, Eskin, Margulis and Mozes [EMM] showed that

$$\frac{\mathrm{vol}_{\mathbb{Z}}\Big( W \cap (R\Omega) \Big)}{\mathrm{vol}\Big( W \cap (R\Omega) \Big)} = 1 + o(1), \qquad \text{as } R \to \infty. \tag{1.14}$$

Furthermore,

$$\mathrm{vol}\big( W \cap (R\Omega) \big) = \lambda (\beta - \alpha) \, R^{d-2} + o(R^{d-2}), \qquad \text{as } R \to \infty,$$

with some $\lambda = \lambda(Q, \Omega) \neq 0$. Eskin, Margulis and Mozes [EMM] provided as well refinements and extensions of (1.14) to lower dimensions.

Introduce the box

$$B(r) = \big\{ x \in \mathbb{R}^d : \, |x|_\infty \leq r \big\}. \tag{1.15}$$

Let $c_0 = c(d, \varepsilon)$ denote a positive constant. Consider the set

$$V \overset{\text{def}}{=} \big\{ Q[x - a] : \ x \in B(r/c_0) \big\} \cap [-c_0 \, r^2, c_0 \, r^2]$$

of values of $Q[x - a]$ lying in the interval $[-c_0 \, r^2, c_0 \, r^2]$, for $x \in B(r/c_0)$. Define the maximal gap between successive values as

$$d(r) \overset{\text{def}}{=} \max_{u \in V} \min \big\{ v - u : \, v > u, \ v \in V \big\}. \tag{1.16}$$

Theorem 1.5. *Let $Q[x]$ be an indefinite quadratic form, $d \geq 9$ and $\varepsilon > 0$. Assume that the constant $c_0 = c(d, \varepsilon)$ is sufficiently small and that $|a| \leq c_0 q^{-1/2} r$. Then the maximal gap satisfies*

$$d(r) \ll_{d,\varepsilon} q^{3d/2} \rho(r^2), \quad \text{for } r^2 \geq c_0^{-1} q^{3d/2},$$

*with $\rho$ defined by* (1.8).

In Section 2 we shall provide as well a bound (see Theorem 2.6) for the remainder term in the quantitative version (1.14) of the Oppenheim conjecture, for $d \geq 9$. This bound is more complicated than the bound of Theorem 1.5 since it depends on the modulus of continuity of the Minkowski functional of the set $\Omega$. In this section we shall mention the following rough Corollaries 1.6 and 1.7 of Theorem 2.6 only.

Corollary 1.6. *Let a quadratic form $Q[x]$ be indefinite and $d \geq 9$. Then, for any $\delta > 0$, there exist* (*sufficiently large*) *constants $C = C(\delta, q, \Omega, d)$ and $C_0 = C_0(\delta, q, \Omega, d)$ such that*

$$(1.17) \quad (1-\delta)\,\mathrm{vol}\big(W \cap (C r \Omega)\big) \leq \mathrm{vol}_{\mathbb{Z}}\big(W \cap (C r \Omega)\big) \leq (1+\delta)\,\mathrm{vol}\big(W \cap (C r \Omega)\big),$$

*provided that*

$$(1.18) \qquad r \geq C_0, \qquad \beta - \alpha \geq C_0, \qquad |a| \leq r, \qquad |\alpha| + |\beta| \leq r^2.$$

Corollary 1.6 is applicable to rational and irrational $Q$. For *irrational $Q$* the approximations can be improved.

Corollary 1.7. *Assume that the quadratic form $Q[x]$ of dimension $d \geq 9$ is irrational and indefinite. Let $R = r T^{1/4}$. Then there exist $T = T(r^2) \to \infty$ such that*

$$(1.19) \quad \Bigg| \frac{\mathrm{vol}_{\mathbb{Z}}\Big(W \cap (R\Omega)\Big)}{\mathrm{vol}\Big(W \cap (R\Omega)\Big)} - 1 \Bigg| \ll_{d,m,q} g(r) + h(r)/(\beta - \alpha) \to 0, \quad \text{as } r \to \infty,$$

*with some functions*

$$g(r) = g(r; q, \Omega, d) \quad \text{and} \quad h(r) = h(r; Q, \Omega, d)$$

*such that $g(r), h(r) \to 0$. The convergence in* (1.19) *holds uniformly in the region where $|a| \leq r$ and $|\alpha| + |\beta| \leq r^2$.*

We have $\rho(s) \to 0$, as $s \to \infty$ (see (1.9)), for irrational $Q$. Thus, Theorem 1.5 gives an upper bound for the maximal gap in Oppenheim's conjecture. The bound of Theorem 1.5 is constructive in the sense that in simple cases one might hope to estimate the quantity $\rho(s)$ explicitly using Diophantine approximation results. In general the estimation of $\rho$ remains an open question. Corollary 1.7 is applicable for shrinking intervals, e.g., for $\beta - \alpha \asymp \sqrt{h(r)}$. The bound of Theorem 2.6 is much more precise than those of Corollaries 1.6 and 1.7. Nevertheless, in order to derive from Theorem 2.6 simple, sharp and precise bounds one needs explicit bounds for $T$, $\gamma(s,T)$ and the modulus of continuity of the functional $M$.

Remark 6.2 shows that the results are uniform over compact sets of irrational matrices $Q$ such that the spectrum of $Q$ is uniformly bounded and is uniformly separated from zero.

The basic steps of the proof consist of:

(1) the introduction of a general approximation problem for the distribution functions of lattice point measures by distribution functions of measures which are absolutely continuous with respect to the Lebesgue measure; both the elliptic as well as hyperbolic cases are obtained as specializations of this general scheme;

(2) an application to the distribution functions of Fourier-Stieltjes transforms, reducing the problem to expansions and integration of Fourier type transforms of the measures (in particular, of certain trigonometric sums) with respect to a one dimensional frequency, say $t$;

(3) integration in $t$ using a basic inequality ([BG1], [BG4]; see (3.12) in this paper), which leads to bounds depending on maximal values $\gamma$ (see (1.7)) of the trigonometric sum;

(4) showing that $\gamma$ tends to zero if and only if the quadratic form is irrational.

Bounds for rates of convergence in the multivariate Central Limit Theorem (CLT) for conic sections (respectively, for bivariate degenerate $U$-statistics) seem to correspond to bounds in the lattice point problems. The "stochastic" diameter (standard deviation) of a sum of $N$ random vectors is of order $\sqrt{N}$, which corresponds to the size of the box of lattice points. In the elliptic case this fact was mentioned by Esseen [Ess], who proved the rate $\mathcal{O}(N^{-1+1/(1+d)})$ for balls around the origin and random vectors with identity covariance, a result similar to the result of Landau [La1]. For sums taking values in a *lattice* and *special* ellipsoids the relation of these error bounds for the lattice point problem and the CLT has been made explicit in Yarnold [Y].

Esseen's result was extended to convex bodies by Matthes [Mat], a result similar to that of Hlawka [H].

The bound $\mathcal{O}(N^{-1})$ in the CLT, for $d \geq 5$, of [BG3] for ellipsoids with diagonal $Q$ and random vectors with independent components (and with arbitrary distribution) is comparable to the results of Jarnik [J1]. The bound $\mathcal{O}(N^{-1})$, for $d \geq 9$, for arbitrary ellipsoids and random vectors — an analogue of the results [BG4] — is obtained in [BG5]. This result is extended to the case of $U$-statistics in [BG6]. Proofs of these probabilistic results are considerably more involved since one has to deal with a more general class of distributions compared to the class of uniform bounded lattice distributions in number theory. A probabilistic counterpart of the results of the present paper remains to be done.

The paper is organized as follows. In Section 2 we formulate the main result, Theorems 2.1 and 2.2, and derive its corollaries and prove the results stated in the introduction. Section 3 is devoted to the proof of Theorems 2.1 and 2.2, using auxiliary results of Sections 4–7. In Section 4 we prove an asymptotic expansion for the Fourier-Stieltjes transforms of the distribution functions and describe some properties of the terms of the expansion. Section 5 contains an integration procedure, which allows to integrate trigonometric sums satisfying the basic inequality (3.12). In Section 6 we obtain a criterion for $Q[x]$ to be irrational in terms of certain trigonometric sums. In Section 7 we investigate the terms of the asymptotic expansions in Theorems 2.1 and 2.2. In Section 8 we obtain auxiliary bounds for the volume of bodies related to indefinite quadratic forms.

We shall use the following notation. By $c$ with or without indices we shall denote generic absolute constants. We shall write $A \ll B$ instead of $A \le cB$. If a constant depends on a parameter, say $d$, then we write $c_d$ or $c(d)$ and use $A \ll_d B$ instead of $A \le c_d B$. By $[B]$ we denote the integer part of a real number $B$.

We shall write $\overline{r} = [r] + 1/2$, for $r \ge 0$. Thus $r \ll \overline{r}$, and for $r \ge 1$ the reverse inequality holds, $\overline{r} \ll r$.

The set of natural numbers is denoted as $\mathbb{N} = \{1, 2, \dots\}$, the set of integer numbers as $\mathbb{Z} = \{0, \pm 1, \pm 2, \dots\}$, and $\mathbb{N}_0 = \{0\} \cup \mathbb{N}$.

We write $B(r) = \big\{ x \in \mathbb{R}^d : \ |x|_\infty \le r \big\}$ and $|x|_\infty = \max\limits_{1 \le j \le d} |x_j|$, $|x|_1 = \sum\limits_{1 \le j \le d} |x_j|$.

The region of integration is specified only in cases when it differs from the whole space. Hence, $\int_{\mathbb{R}} = \int$ and $\int_{\mathbb{R}^d} = \int$.

We use the notation

$$\mathrm{e}\big\{ t \big\} = \exp\big\{ it \big\}, \qquad i = \sqrt{-1}, \tag{1.20}$$

which differs by an inessential factor $2\pi$ from often used $\mathrm{e}\big\{ t \big\} = \exp\big\{ 2\pi i t \big\}$. Since we study forms with arbitrary real coefficients, the convention (1.20) suppresses lots of immaterial factors $2\pi$.

The Fourier-Stieltjes transforms of functions, say $F : \mathbb{R} \to \mathbb{R}$, of bounded variation are denoted as

$$\widehat{F}(t) = \int \mathrm{e}\big\{ ts \big\} \, dF(x).$$

Throughout $\mathbf{I}\{A\}$ denotes the indicator function of event $A$, that is, $\mathbf{I}\{A\} = 1$ if $A$ occurs, and $\mathbf{I}\{A\} = 0$ otherwise.

For $s > 0$, define the function

$$\mathcal{M}(t;s) = \big(|t|\,s\big)^{-1}\,\mathbf{I}\big\{\,|t| \le s^{-1/2}\big\} + |t|\,\mathbf{I}\big\{\,|t| > s^{-1/2}\big\}. \tag{1.21}$$

For a multi-index $\alpha = (\alpha_1, \dots, \alpha_d)$, we write $\alpha! = \alpha_1! \dots \alpha_d!$. Partial derivatives of functions $f : \mathbb{R}^d \to \mathbb{C}$ we denote by

$$\partial^\alpha f(x) = \partial_x^\alpha f(x) = \frac{\partial^{\alpha_1}}{(\partial x_1)^{\alpha_1}} \cdots \frac{\partial^{\alpha_d}}{(\partial x_d)^{\alpha_d}} f(x).$$

Sometimes we shall use notation related to Fréchet derivatives: for $\alpha = (\alpha_1, \dots, \alpha_n)$, we write

$$f^{(|\alpha|_1)}(x) h_1^{\alpha_1} \dots h_n^{\alpha_n} = \partial_{t_1}^{\alpha_1} \dots \partial_{t_n}^{\alpha_n} f(x + t_1 h_1 + \dots + t_n h_n)\Big|_{t_1 = \dots = t_n = 0}. \tag{1.22}$$

*Acknowledgment.* We would like to thank G. Margulis for drawing our attention to the close relation between the quantitative Oppenheim conjecture and the lattice point remainder problem and helpful discussions. Furthermore, we would like to thank A.Yu. Zaitsev for a careful reading of the manuscript and useful comments.

## 2. The main result: Proofs of the theorems of the introduction

For the formulation of the main result, Theorem 2.1, we need some simple notions related to measures on $\mathbb{R}^d$. We shall consider *signed* measures, that is, $\sigma$-additive set functions $\mu : \mathcal{B}^d \to \mathbb{R}$, where $\mathcal{B}^d$ denotes the $\sigma$-algebra of Borel subsets of $\mathbb{R}^d$. Probability measure (or distribution) is a nonnegative and normalized measure (that is, $\mu(C) \ge 0$, for $C \in \mathcal{B}^d$ and $\mu(\mathbb{R}^d) = 1$). We shall write $\int f(x)\,\mu(dx)$ for the (Lebesgue) integral over $\mathbb{R}^d$ of a measurable function $f : \mathbb{R}^d \to \mathbb{C}$ with respect to a signed measure $\mu$, and denote as usual by $\mu * \nu(C) = \int \mu(C - x)\,\nu(dx)$, for $C \in \mathcal{B}^d$, the convolution of the signed measures $\mu$ and $\nu$. Equivalently, $\mu * \nu$ is defined as the signed measure such that

$$\int f(x)\,\mu * \nu(dx) = \iint f(x+y)\,\mu(dx)\,\nu(dy), \tag{2.1}$$

for any integrable function $f$.

Let $p_x \in \mathbb{R}$, $x \in \mathbb{Z}^d$, be a system of weights. Using signed measures, weighted trigonometric sums, say,

$$\sum_{x \in \mathbb{Z}^d} \mathrm{e}\big\{\,t\,Q[x]\,\big\}\,p_x = \int \mathrm{e}\big\{\,t\,Q[x]\,\big\}\,\theta(dx), \qquad \mathrm{e}\{v\} = \exp\{iv\},$$

can be represented as an integral with respect to the signed measure $\theta$ concentrated on the lattice $\mathbb{Z}^d$ such that $\theta\big(\{x\}\big) = p_x$, for $x \in \mathbb{Z}^d$.

The *uniform lattice measure* $\mu(\,\cdot\,;r)$ concentrated on the lattice points in the cube $B(r)=\big\{\,x\in\mathbb{R}^d:\ |x|_\infty\le r\,\big\}$ is defined by

$$\mu(C;r)=\frac{\mathrm{vol}_{\mathbb{Z}}\Big(C\cap B(r)\Big)}{\mathrm{vol}_{\mathbb{Z}}B(r)},\qquad \text{for}\quad C\in\mathcal{B}^d. \tag{2.2}$$

In other words, the measure $\mu(\,\cdot\,;r)$ assigns equal weights $\mu(\{x\};r)=(2\overline{r})^{-d}$ to lattice points in the cube $B(r)$, where $\overline{r}=[r]+1/2$. Notice as well that $\mu(\,\cdot\,;r)=\mu(\,\cdot\,;\overline{r})$.

We define the *uniform measure* $\nu(\,\cdot\,;r)$ in $B(r)$ by

$$\nu(C;r)=\frac{\mathrm{vol}\Big(C\cap B(r)\Big)}{\mathrm{vol}\,B(r)},\qquad \text{for}\quad C\in\mathcal{B}^d. \tag{2.3}$$

For a number $R>0$ write $\Phi=\mu(\,\cdot\,;R)$ and $\Psi=\nu(\,\cdot\,;\overline{R})$, and introduce the measures

$$\mu=\Phi*\mu^{*k}(\,\cdot\,;r),\quad \nu=\Psi*\nu^{*k}(\,\cdot\,;\overline{r}),\quad k\in\mathbb{N}. \tag{2.4}$$

The distribution function, say $G$, of a quadratic form $Q[x-a]$ with respect to a signed measure, say $\lambda$, on $\mathbb{R}^d$ is defined as

$$G(s)=\lambda\big\{\,x\in\mathbb{R}^d:\ Q[x-a]\le s\,\big\}=\int\mathbf{I}\big\{\,Q[x-a]\le s\,\big\}\,\lambda(dx), \tag{2.5}$$

where $\mathbf{I}\big\{\,A\,\big\}$ denotes the indicator function of event $A$. The function $G:\mathbb{R}\to\mathbb{R}$ is right continuous and satisfies $G(-\infty)=0$, $G(\infty)=\lambda(\mathbb{R}^d)$. If $\lambda$ is a probability measure (i.e., nonnegative and normalized) then we have in addition: $G:\mathbb{R}\to[0,1]$ is nondecreasing and $G(\infty)=1$.

We shall obtain an asymptotic expansion of the distribution function, say $F$, of $Q[x-a]$ with respect to the measure $\mu$ defined by (2.4). The first term of this expansion will be the distribution function, say $F_0$, of $Q[x-a]$ with respect to the measure $\nu$ defined by (2.4). Other terms of this expansion will be distribution functions $F_j$, $j\in2\mathbb{N}$, of certain signed measures related to the measure $\nu$ (or, in other words, to certain Lebesgue type volumes). A description of $F_j$ will be given after Theorem 2.2.

Introduce the function (cf. (1.4))

$$\varphi_a(t;r^2)=\Big|\int\mathrm{e}\big\{\,tQ[x-a]\,\big\}\,\mu^{*3}(dx;r)\Big|, \tag{2.6}$$

and, for a number $T\ge1$, define (cf. (1.6) and (1.7))

$$\gamma\big(r^2,T\big)\overset{\mathrm{def}}{=}\sup_{a\in\mathbb{R}^d}\ \sup_{r^{-1}\le t\le T}\ \varphi_a(t;r^2). \tag{2.7}$$

Our main result is the following theorem.

THEOREM 2.1. *Assume that*

$$d \geq 9, \quad p \in \mathbb{N}, \quad 2 \leq p < d/2, \quad k \geq 2p+2, \quad 0 \leq r \leq R, \quad T \geq 1.$$

*Then the distribution function $F$ allows the following asymptotic expansion*

$$F(s) = F_0(s) + \sum_{j \in 2\mathbb{N},\ j<p} F_j(s) + \mathcal{R} \tag{2.8}$$

*with a remainder term $\mathcal{R}$ satisfying*

$$|\mathcal{R}| \ll_{d,k,\varepsilon} \frac{q^{d/2}}{r^2 T} + \frac{R^p}{r^{2p}} \Big(1 + \frac{|a|}{r}\Big)^p q^{p+d/2} + \gamma^{1-8/d-\varepsilon}\big(r^2, T\big)\, T^{\varepsilon}\, \frac{q^{d/2}}{r^2}, \tag{2.9}$$

*for any $\varepsilon > 0$.*

Notice, that the estimate (2.9) is uniform in $s$.

The measure $\Phi$ (or its support $B(R) \cap \mathbb{Z}^d$) represents the main box of size $R$ from which lattice points are taken. The convolution of $\Phi$ with $\mu^{*k}(\cdot\,; r)$ is a somewhat smoother lattice measure than $\Phi$. Note though that the weights assigned by $\mu$ to the lattice points near the boundary of the box $B(R)$ become smaller when the points approach the boundary of the box $B(R+kr)$. The weights assigned to lattice points in $B(R-kr)$ remain unchanged. Later on we shall choose the size $r$ of the smoothing measure $\mu(\cdot\,; r)$ smaller in comparison with $R$, that is, we shall assume that $R \geq ckr$ with a sufficiently large constant $c$. This smoothing near the boundary simplifies the derivation of approximations and helps to avoid extra logarithmic factors in the estimates of errors. The corresponding measure $\nu$ is the continuous counterpart of $\mu$ with the dominating counting measure on $\mathbb{Z}^d$ replaced by the Lebesgue measure on $\mathbb{R}^d$. Theorem 2.1 allows a generalization. The measure $\Phi$ can be replaced by an arbitrary uniform lattice measure with support in a cube of size $R$, see Theorem 2.2 below. Theorem 2.1 is a partial case of Theorem 2.2. We shall prove Theorem 2.2 in Section 3. In order to formulate that result, we extend our notation.

We shall denote $\pi = \nu(\cdot\,; 1/2)$. The measure $\pi$ has the density $\frac{d\pi}{dx} = \mathbf{I}\big\{ |x|_\infty \leq 1/2 \big\}$ with respect to the Lebesgue measure in $\mathbb{R}^d$, so that $\pi(dx) = \mathbf{I}\big\{ |x|_\infty \leq 1/2 \big\}\, dx$. The measure $\nu(\cdot\,; \overline{r})$ has the density $(2\overline{r})^{-d}\, \mathbf{I}\big\{ |x|_\infty \leq \overline{r} \big\}$. Notice as well that $\nu(\cdot\,; \overline{r}) = \pi * \mu(\cdot\,; \overline{r})$.

Henceforth $\Phi$ will denote a probability measure on $\mathbb{R}^d$ such that $\Phi(A) = 1$, for some subset $A \subset B(R) \cap \mathbb{Z}^d$ and

$$\Phi\big(\{x\}\big) = 1/\operatorname{card} A, \qquad \text{for all } x \in A. \tag{2.10}$$

We do not impose restrictions on the structure of $A$ except that $A \subset B(R) \cap \mathbb{Z}^d$ and $A \neq \emptyset$. Write $\Psi = \Phi * \pi$. It is easy to see that $\Psi$ has the density

$$\frac{d\Psi}{dx} = \frac{1}{\operatorname{card} A} \sum_{y \in A} \mathbf{I}\big\{ |x-y|_\infty \leq 1/2 \big\}. \tag{2.11}$$

We define measures $\mu$ and $\nu$ as in (2.4), and denote distribution functions of $Q[x-a]$ with respect to $\mu$ and $\nu$ as $F$ and $F_0$ respectively. Notice that the measure $\nu$ has the density

$$D(x) \stackrel{\text{def}}{=} \frac{d\nu}{dx} = \frac{d\Psi}{dx} * \Big(\frac{d\nu(\,\cdot\,;\overline{r})}{dx}\Big)^{*k}, \tag{2.12}$$

where $f * g$ denotes the convolution of functions $f$ and $g$,

$$f * g(x) = \int f(x-y)\, g(y)\, dy.$$

Using the Fourier transform, we can easily verify that the density $D$ admits continuous bounded partial derivatives $|\partial^\alpha D(x)| \ll_{d,k} \overline{r}^{\,-d-|\alpha|_1}$, for $|\alpha|_\infty \le k-2$ (see Lemma 7.1).

Theorem 2.2. *Theorem* 2.1 *holds with* $\Phi$ *and* $\Psi$ *defined by* (2.10) *and* (2.11) *respectively.*

Let us now define the functions $F_j$, for $j \in 2\mathbb{N}$. Let $\eta = (\eta_1, \dots, \eta_m)$ denote a multi-index with entries $\eta_1, \dots, \eta_m \in \mathbb{N}$. Write $\sum^{**}_{\eta\colon\, |\eta|_1 = j}$ for the sum which extends over all possible representations of the *even* number $j$ as a sum $j = \eta_1 + \dots + \eta_m$ of *even* $\eta_1, \dots, \eta_m \ge 2$, for all possible $m \ge 1$. For example, for $j = 6$, we have $6 = 6$, $6 = 4+2$, $6 = 2+4$ and $6 = 2+2+2$. Introduce the functions

$$D_j(x) = \sum^{**}_{\eta\colon\, |\eta|_1 = j} D_{j\eta}(x) \tag{2.13}$$

with

$$D_{j\eta}(x) = \frac{(-1)^m}{\eta!} \int \cdots \int D^{(j)}(x) u_1^{\eta_1} \dots u_m^{\eta_m} \prod_{l=1}^{m} \pi^{*(k+1)}(du_l), \tag{2.14}$$

where the density $D$ is defined by (2.12), and where we use the notation (1.22) for the Fréchet derivatives. For example, we have

$$D_2(x) = -\tfrac{1}{2} \int D''(x) u^2\, \pi^{*(k+1)}(du),$$

and $D_4(x) = D_{44}(x) + D_{422}(x)$ with

$$D_{44}(x) = -\tfrac{1}{24} \int D^{(4)}(x) u^4\, \pi^{*(k+1)}(du),$$
$$D_{422}(x) = \tfrac{1}{4} \iint D^{(4)}(x) u_1^2 u_2^2\, \pi^{*(k+1)}(du_1)\, \pi^{*(k+1)}(du_2).$$

Let $\nu_j$ denote the signed measure on $\mathbb{R}^d$ with density $D_j$. We define the function $F_j$, for $j \in 2\mathbb{N}$, as the distribution function of $Q[x-a]$ with respect to the signed measure $\nu_j$; that is,

$$F_j(s) = \nu_j\big\{x \in \mathbb{R}^d :\ Q[x-a] \le s\big\} = \int \mathbf{I}\big\{Q[x-a] \le s\big\}\, D_j(x)\, dx. \tag{2.15}$$

The function $F_j : \mathbb{R} \to \mathbb{R}$ is a function of bounded variation, $F_j(-\infty) = F_j(\infty) = 0$ and

$$\sup_s \bigl| F_j(s) \bigr| \ll_{j,d} \tfrac{R^j}{r^{2j}} \Bigl(1 + \tfrac{|a|}{r}\Bigr)^j q^{j+d/2}, \qquad \text{for } j < d/2; \tag{2.16}$$

see Lemma 7.4.

In the elliptic case the choice of $\Phi$ is immaterial as long as the support of $\Phi$ contains a sufficiently massive box of lattice points. Thus we shall simply choose $\Phi$ and $\Psi$ as in (2.4). The same choice of $\Phi$ is appropriate for the estimation of maximal gaps (cf. Theorem 1.5) in the hyperbolic case. The choice of a general as possible $\Phi$ is appropriate for proving refinements of (1.14). We shall restrict ourselves to the following special $\Phi$ generated by a star-shaped closed bounded set $\Omega$ (see (1.13)) whose nonempty interior contains zero. Define

$$\Phi(C) = \frac{\operatorname{vol}_{\mathbb{Z}}\Bigl( C \cap (R\Omega) \Bigr)}{\operatorname{vol}_{\mathbb{Z}}(R\Omega)} \tag{2.17}$$

and let in accordance with (2.10) the set $A$ be given by $A = (R\,\Omega) \cap \mathbb{Z}^d$. The measure $\Psi$ is again defined by (2.11). In order to guarantee that $\Bigl\{ x \in \mathbb{Z}^d : \Phi\bigl(\{x\}\bigr) > 0 \Bigr\} \subset B(R)$, we shall assume throughout that $\Omega \subset B(1)$; this is not a restriction of generality. Hence, for the Minkowski functional of the set $\Omega$ we have

$$|x|_\infty \le \bigl| M(x) \bigr| \le m\, |x|_\infty, \qquad \text{for all } x \in \mathbb{R}^d, \tag{2.18}$$

with some $m \ge 1$. The inequalities (2.18) are equivalent to $B(1/m) \subset \Omega \subset B(1)$.

The modulus of continuity

$$\omega(\delta) = \sup_{|y|_\infty \le \delta,\ |x|_\infty = 1} \bigl| M(x+y) - M(x) \bigr| \tag{2.19}$$

of $M$ satisfies $\lim_{\delta \to 0} \omega(\delta) = 0$. For $\Phi$ and $\Psi$ in (2.4) we have $\Omega = B(1)$, $M(x) = |x|_\infty$ and $\omega(\delta) = \delta$.

Let $(\partial\Omega)_\sigma \overset{\text{def}}{=} \partial\Omega + B(\sigma)$ be a $\sigma$-neighborhood of the boundary $\partial\Omega$ of $\Omega$. Then, introducing the weight $p_0 = 1/\operatorname{vol}_{\mathbb{Z}}(R\Omega)$, writing for a while $\sigma = kr/R$ and assuming that $\Phi$ is defined by (2.17), we have

$$\begin{aligned} &\mu(C) = p_0 \operatorname{vol}_{\mathbb{Z}} C, && \text{for } C \subset R\bigl(\Omega \setminus (\partial\Omega)_\sigma\bigr), \\ &0 \le \mu(C) \le p_0 \operatorname{vol}_{\mathbb{Z}}\bigl(C \cap (R\,\Omega_\sigma)\bigr), && \text{for } C \subset \mathbb{R}^d, \\ &\mu(C) = 0, && \text{for } C \subset \mathbb{R}^d \setminus (R\,\Omega_\sigma). \end{aligned} \tag{2.20}$$

Notice that $p_0 = (2\overline{R})^{-d}$ for $\Phi$ defined after (2.4). In order to prove (2.20), it suffices to consider the case when the set $C$ is a one point set, and to use elementary properties of convolutions. Similarly, for measurable $C \subset \mathbb{R}^d$, we have

$$
\begin{aligned}
&\nu(C) = p_0 \int_C dx, && \text{for } C \subset R\big(\Omega \setminus (\partial\Omega)_\sigma\big), \\
&0 \le \nu(C) \le p_0 \int_{C\cap(R\,\Omega_\sigma)} dx, && \text{for } C \subset \mathbb{R}^d, \\
&\nu(C) = 0, && \text{for } C \subset \mathbb{R}^d \setminus (R\Omega_\sigma),
\end{aligned} \tag{2.21}
$$

where now $\sigma = (kr+1)/R$. Using (2.19), it is easy to see that

$$
R(\partial\Omega)_\sigma \subset \big\{x \in \mathbb{R}^d :\ 1 - \omega(\sigma) \le M(x/R) \le 1 + \omega(\sigma)\big\}, \quad \text{for any } \sigma > 0. \tag{2.22}
$$

We conclude the section by deriving all results of the introduction as corollaries of Theorem 2.1; the refinement of (1.14), Theorem 2.6, is implied by Theorem 2.2.

*The elliptic case.* Let us start with the following corollary of Theorem 2.1.

COROLLARY 2.3. *Assume that the operator $Q$ is positive and $T \ge 1$. Then we have*

$$
\big|\operatorname{vol}_{\mathbb{Z}}(E_s+a) - \operatorname{vol} E_s\big| \ll_{d,\varepsilon} (s+1)^{d/2}\, q^{p+d/2} \Big(\frac{1}{s^{p/2}} + \frac{1}{sT} + \gamma^{1-8/d-\varepsilon}\big(s,T\big)\,\frac{T^\varepsilon}{s}\Big), \tag{2.23}
$$

*for $d \ge 9$, $2 \le p < d/2$, $p \in \mathbb{N}$ and any $\varepsilon > 0$. The quantity $\gamma(s,T)$ is defined in* (1.7) (*cf.* (2.6) *and* (2.7)).

*Proof.* The bound (2.23) obviously holds for $s \le 1$. Therefore proving (2.23) we shall assume that $s > 1$.

We assume as well that $|a|_\infty \le 1$. This assumption does not restrict generality since

$$
\operatorname{vol}_{\mathbb{Z}}(E_s + a) = \operatorname{vol}_{\mathbb{Z}}(E_s + a - m), \qquad \text{for any } m \in \mathbb{Z}^d, \tag{2.24}
$$

and in (2.23) we can replace $a$ by $a - m$ with some $m \in \mathbb{Z}^d$ such that $|a - m|_\infty \le 1$.

Choose the measure $\Phi$ as in (2.4), and

$$
k = 2p + 2, \qquad r = \sqrt{s}, \qquad R = 2kr.
$$

Obviously $\overline{R} \ll_k \sqrt{s+1}$, for $s > 1$. Therefore the bound (2.9) of Theorem 2.1 implies (2.23) provided that we verify that our choices yield

$$
F(s) = (2\,\overline{R})^{-d} \operatorname{vol}_{\mathbb{Z}}(E_s + a), \quad F_0(s) = (2\,\overline{R})^{-d} \operatorname{vol} E_s, \quad F_j(s) = 0, \ \text{for } j \ne 0. \tag{2.25}
$$

Let us prove the first equality in (2.25). The ellipsoid $E_1$ is contained in the unit ball, that is, $E_1 \subset \{|x| \le 1\} \subset B(1)$, since the modulus of the minimal eigenvalue $q_0$ is 1. Therefore we have $E_s \subset B(\sqrt{s})$. Due to our choice of $R$, $r$ and $k \ge 6$, we have $R - kr \ge 6\sqrt{s}$. Thus, the inequality $|a|_\infty \le 1$, the relations $E_s + a \subset B(1+\sqrt{s})$ and $F(s) = \mu(E_s + a)$ together with (2.20) imply the first equality in (2.25).

Let us prove the second equality in (2.25). Using (2.12) and (2.21) we see that the density $D$ is equal to zero outside the set $B(R + kr + 1)$, and $D(x) \equiv (2\overline{R})^{-d}$, for $x \in B(R - kr - 1)$. The ellipsoid $E_s + a$ is a subset of $B(R - kr - 1)$, that yields the second equality in (2.25).

For the proof of $F_j(s) = 0$ notice that the derivatives of $D$ vanish in $B(R - kr - 1)$, hence in the ellipsoid $E_s + a$ as well. □

*Proof of Theorem* 1.3. This theorem is implied by Corollary 2.3. Indeed, the estimate (1.10) is obvious for $s \le 1$. For $s > 1$, the estimate (1.10) is implied by (2.23) estimating $q^p \le q^{d/2}$, choosing $p = 2\zeta$ and $T = T(s)$ as in the condition of Theorem 1.3. □

*Proof of Corollary* 1.4. We have to prove (1.11) and (1.12). The proof of (1.12) reduces to proving that

$$\operatorname{vol}_{\mathbb{Z}}(E_{\tau+\delta} + a) - \operatorname{vol}_{\mathbb{Z}}(E_\tau + a) > 0,$$

for $\delta \ge c(d,\varepsilon)\, q^{3d/2}\, \rho_0(s)$ with a sufficiently large constant $c(d,\varepsilon)$. Using (1.11) it suffices to verify that $|\mathcal{R}| \le 1/2$, which is obviously fulfilled.

Let us prove (1.11). Consider an interval $(\tau, \tau+\delta]$ with $\tau \ge s \ge 1$. We shall apply the bound of Theorem 1.3 which for $s \ge 1$ yields

$$\Big|\operatorname{vol}_{\mathbb{Z}}(E_s + a) - \operatorname{vol} E_s\Big| \ll_{d,\varepsilon} q^d\, s^{-1+d/2}\, \rho_0(s). \tag{2.26}$$

We get

$$\Big|\operatorname{vol}_{\mathbb{Z}}\big((E_{\tau+\delta} + a) \setminus (E_\tau + a)\big) - \operatorname{vol}\big(E_{\tau+\delta} \setminus E_\tau\big)\Big| \ll_{d,\varepsilon} q^d\, (\tau+\delta)^{d/2-1}\, (\rho_0(\tau+\delta) + \rho_0(\tau)). \tag{2.27}$$

The estimate (2.27) implies (1.11). Just note that $\rho_0(\tau) \le \rho_0(s)$, for $\tau \ge s$, divide both sides of (2.27) by $\operatorname{vol}\big(E_{\tau+\delta} \setminus E_\tau\big)$ and use

$$\operatorname{vol}\big(E_{\tau+\delta} \setminus E_\tau\big) = \big((\tau+\delta)^{d/2} - \tau^{d/2}\big) \operatorname{vol} E_1,$$

$$(\tau+\delta)^{d/2} - \tau^{d/2} \gg_d \int_{\tau+\delta/2}^{\tau+\delta} u^{-1+d/2}\, du \gg_d \delta\, (\tau+\delta/2)^{-1+d/2} \gg_d \delta\, (\tau+\delta)^{-1+d/2},$$

$$\operatorname{vol} E_1 = \operatorname{vol}\big\{x \in \mathbb{R}^d : Q[x] \le 1\big\} \ge \operatorname{vol}\big\{x \in \mathbb{R}^d : |x| \le 1/\sqrt{q}\big\} = c_d\, q^{-d/2}.$$

□

*Proof of Corollary* 1.2. It suffices to use (1.9) and (1.12). □

*Proof of Theorem* 1.1. Assuming the irrationality of $Q$, the bound $o(s^{-1})$ is implied by Theorem 1.3 and (1.9) since $\operatorname{vol} E_s \gg_d q^{-d/2}\, s^{d/2}$.

The bound $o(s^{-1})$ in (1.2) implies $d(\tau, Q, 0) \to 0$, as $\tau \to \infty$, which is impossible for rational $Q$. □

*The hyperbolic case.* For an interval $I = (\alpha, \beta] \subset \mathbb{R}$ we write

$$F(I) = F(\beta) - F(\alpha), \qquad F_j(I) = F_j(\beta) - F_j(\alpha). \tag{2.28}$$

Notice that $F(I) = F(s)$ in the case $I = (-\infty, s]$. Theorem 2.2 has the following obvious corollary.

Corollary 2.4. *Under the conditions of Theorem* 2.2 *we have*

$$F(I) = F_0(I) + \sum_{j \in 2\mathbb{N},\ j<p} F_j(I) + \mathcal{R} \tag{2.29}$$

*with remainder term* $\mathcal{R}$ *which satisfies* (2.9).

Lemma 2.5. *Let* $Q[x]$ *be an indefinite quadratic form of dimension* $d \geq 9$. *Let* $p$, $k$ *and* $\varepsilon$ *satisfy the conditions of Theorem* 2.1. *Assume that* $M(x) = |x|_\infty$ *and* $\Omega = B(1)$. *Let* $c_1 = c_1(d, \varepsilon)$ *denote a sufficiently small positive constant. Finally, assume that*

$$R \geq r \geq \tfrac{1}{c_1}, \quad \alpha, \beta \in [-c_1 R^2, c_1 R^2], \quad \tfrac{r}{R} \leq \tfrac{c_1}{k}, \quad \sqrt{q}\, |a| \leq c_1 R. \tag{2.30}$$

*Then*

$$\Big| \frac{F(I)}{F_0(I)} - 1 \Big| \ll_{p,k,d,\varepsilon} \frac{q^{3d/2}}{r^2}\, \frac{R^d}{r^d} + \frac{q^{p+d} R^2}{\beta - \alpha} \Big( \frac{1}{r^2 T} + \frac{R^{2p}}{r^{3p}} + \gamma^{1-8/d-\varepsilon}\big(r^2, T\big)\, \frac{T^\varepsilon}{r^2} \Big). \tag{2.31}$$

*Proof.* The result follows from Corollary 2.4 dividing (2.29) by $F_0(I)$ and using the estimates

$$\sum_{j \in 2\mathbb{N},\ j<p} \big| F_j(I) \big| \ll_{k,d} (\beta - \alpha)\, q^{d-2}\, r^{-2-d}\, R^{d-2}, \tag{2.32}$$

$$F_0(I) \gg_{k,d} (\beta - \alpha)\, q^{-d/2}\, R^{-2}. \tag{2.33}$$

Let us prove (2.32). Using (2.15), (2.28), (2.30), (7.1), applying the estimate (8.9) of Lemma 8.2 with

$$M(x) = |x|_\infty, \quad R = \overline{R} + k\,\overline{r}, \quad \text{and} \quad m = \lambda = 1,$$

using the bound $|a_0| \le \sqrt{q}\,|a|$, we obtain

$$
\begin{aligned}
\tag{2.34}
\bigl|F_j(I)\bigr| &\ll_{k,d} \overline{r}^{\,-j-d} \int \mathbf{I}\bigl\{ |x|_\infty \le \overline{R} + k\overline{r} \bigr\}\, \mathbf{I}\bigl\{ Q[x-a] \in I \bigr\}\, dx \\
&\ll_{k,d} (\beta-\alpha)\, q^{(d-2)/2}\, \overline{r}^{\,-j-d} \Bigl(1 + \tfrac{\sqrt{q}\,|a|}{\overline{R}+k\,\overline{r}}\Bigr)^{d-2} (\overline{R}+k\,\overline{r})^{d-2} \\
&\ll_{k,d} (\beta-\alpha)\, q^{d-2}\, r^{-j-d}\, R^{d-2}.
\end{aligned}
$$

In the proof of (2.34) the condition $R \ge r \ge 1/c_1$ allowed us to replace $\overline{R}$ and $\overline{r}$ by $R$ and $r$ respectively. Summation in $j$, $2 \le j < p$, of the inequalities (2.34) yields (2.32).

Let us prove (2.33). Using (2.12), (2.21), the lower bound (8.10) of Lemma 8.2 and conditions (2.30), we have

$$
\begin{aligned}
F_0(I) \gg_d R^{-d} \int \mathbf{I}\bigl\{ |x|_\infty \\
\le \overline{R} - k\overline{r} \bigr\}\, \mathbf{I}\bigl\{ Q[x-a] \in I \bigr\}\, dx \gg_d (\beta-\alpha)\, q^{-d/2}\, R^{-2}.
\end{aligned}
$$

□

*Proof of Theorem* 1.5. We shall apply Lemma 2.5 choosing $T = T(r^2)$, $R = r\,k/c_1$, the maximal $p$ and minimal $k$ such that the conditions of Theorem 2.1 are fulfilled. In this particular case we can rewrite (2.31) as

$$
\tag{2.35}
\Bigl| \frac{F(I)}{F_0(I)} - 1 \Bigr| \ll_{d,\varepsilon} q^{3d/2}\, r^{-2} + \frac{q^{d+p}}{\beta-\alpha}\, \rho(r^2).
$$

In order to estimate the maximal gap, it suffices to show that any interval $I = (\alpha, \beta]$ contains at least one value of the quadratic form (i.e., $F(I) > 0$) provided that $\alpha$ and $\beta$ satisfy $\beta - \alpha \gg_{d,\varepsilon} q^{d+p}\, \rho(r^2)$. The inequality $F(I) > 0$ holds if the right-hand side of (2.35) is bounded from above by a sufficiently small constant which can depend on $d$ and $\varepsilon$. □

Next we formulate and prove some refinements of (1.14). Recall that the set $\Omega$ satisfies $B(1/m) \subset \Omega \subset B(1)$ (see (2.18)), and that $\omega$ denotes the modulus of continuity of the Minkowski functional $M$ of the set $\Omega$ (see (1.13) and (2.19)). Write

$$
\tag{2.36}
W = \bigl\{ x \in \mathbb{R}^d :\ Q[x-a] \in I \bigr\}, \qquad I = (\alpha, \beta],
$$

$$
\tag{2.37}
\Delta^* = \Biggl| \frac{\mathrm{vol}_{\mathbb{Z}}\Bigl( W \cap (R\Omega) \Bigr)}{\mathrm{vol}\Bigl( W \cap (R\Omega) \Bigr)} - 1 \Biggr|,
$$

and

$$
\tag{2.38}
\begin{aligned}
&\varepsilon_0 = r/R, \qquad && \varepsilon_1 = \sqrt{q}\,|a|/R, \qquad && \varepsilon_2 = (|\alpha| + |\beta|)/R^2, \\
&\varepsilon_3 = \omega(\varepsilon_0/c_2), && \varepsilon_4 = \omega(\varepsilon_1\sqrt{q}/c_2), && \varepsilon_5 = \omega(\varepsilon_2\sqrt{q}/c_2),
\end{aligned}
$$

where $c_2 = c_2(d, m)$ denotes a positive constant.

Theorem 2.6. *Assume that the form* $Q[x]$ *is indefinite and* $d \geq 9$. *Let* $r \geq 1/c_2$, $0 < \varepsilon < 1 - 8/d$, $T \geq 1$ *and*

$$\varepsilon_j \leq c_2, \ \text{for } j = 0,1,2,3, \quad \varepsilon_j \leq c_2\, q^{-1/2}, \ \text{for } j = 4,5. \tag{2.39}$$

*Then we have*

$$\Delta^* \ll_{d,\varepsilon,m} q^{d-1}\Big(1 + \frac{1}{r^2\,\varepsilon_0^d}\Big)(\varepsilon_3 + \varepsilon_4 + \varepsilon_5) + \frac{q^{3d/2}}{\beta - \alpha}\Big(\frac{1}{\varepsilon_0^2 T} + \frac{1}{r^{-3+d/2}\,\varepsilon_0^{d+2}} + \gamma^{1-8/d-\varepsilon}\big(r^2,T\big)\,\frac{T^\varepsilon}{\varepsilon_0^2}\Big) \tag{2.40}$$

*provided that the constant* $c_2$ *is sufficiently small.*

*Proof of Corollary* 1.6. Let us apply Theorem 2.6 choosing $c_2 = c_2(\delta, q, d, m)$ sufficiently small depending on $\delta$ and $q$ as well. Choose $R = C\,r$. Then we have $\varepsilon_0 = C^{-1}$ and the last two inequalities in (1.18) guarantee conditions (2.39). In particular, we have $\varepsilon_j \leq c_2$, for $j = 3,4,5$. Hence, we can apply the bound (2.40). Choosing $T = 1$, $\varepsilon = (1 - 8/d)/2$ and estimating $\gamma \leq 1$, we obtain

$$\Delta^* \ll_{q,d,m} c_2\,(1 + r^{-2} C^d) + (C^2 + r^{-1}\, C^{d+2})/(\beta - \alpha), \tag{2.41}$$

for $r \geq 1$. The estimates (1.17) follow if the right-hand side of (2.41) is bounded from above by $\delta$. But this holds in view of the first two inequalities in (1.18) and our choices of $C$, $C_0$ and $c_2$. □

*Proof of Corollary* 1.7. During the proof we shall write $T = T(r^2)$ and $\gamma = \gamma\big(r^2, T\big)$. Choose $\varepsilon = (1 - 8/d)/2$. Since $R = r T^{1/4}$ with $T \to \infty$, we have $\varepsilon_0 = T^{-1/4}$ and the conditions (2.39) are fulfilled for sufficiently large $r$. Moreover, $\varepsilon_j \to 0$, for $j = 3,4,5$. Hence, the bound (2.40) yields (1.19) with

$$g(r) = (1 + r^{-2}\, T^{d/4})\,(\varepsilon_3 + \varepsilon_4 + \varepsilon_5),$$
$$h(r) = T^{-1/2} + r^{-1}\, T^{(d+2)/4} + \gamma^{(1-8/d)/2}\, T^{1-4/d}.$$

If $g(r) \nrightarrow 0$ or $h(r) \nrightarrow 0$, we can choose (if necessary) $T = T(r^2) \to \infty$ growing somewhat slower such that $g(r) \to 0$ or $h(r) \to 0$ (cf. a similar redefinition of $T(s)$ in the case of (1.9)). □

For the proof of Theorem 2.6 we shall need the following Lemma 2.7, which allows us to estimate $\Delta^*$ using Theorem 2.2. Write

$$\Delta^\circ = \mathrm{vol}_{\mathbb{Z}}\big(W \cap (R\Omega)\big) - \mathrm{vol}\big(W \cap (R\Omega)\big). \tag{2.42}$$

By $\mu_R = \Phi * \mu^{*k}(\cdot\,; r)$ and $\nu_R = \Psi * \mu^{*k}(\cdot\,; r)$, where $\Psi = \Phi * \pi$, we shall denote the measures $\mu$ and $\nu$, emphasizing the dependence on the parameter $R$ which enters in the definition (2.17) of $\Phi$.

LEMMA 2.7. *Let* $\sigma_0 = (kr+1)/R$. *Assume that* $\omega(2\sigma_0) < 1$. *Write*

$$s_1 = R/\bigl(1+\omega(2\sigma_0)\bigr), \quad s_2 = R/\bigl(1-\omega(2\sigma_0)\bigr), \quad p_2 = 1/\operatorname{vol}_{\mathbb{Z}}(s_2\Omega). \tag{2.43}$$

*Then we have*

$$|\Delta^\circ| \le p_2^{-1} \max_{s=s_1,s_2} |\Delta_s| + v, \tag{2.44}$$

*where*

$$v = \operatorname{vol}\Bigl( W \cap \bigl\{ s_1\bigl(1-\omega(2\sigma_0)\bigr) \le M(x) \le s_2\bigl(1+\omega(2\sigma_0)\bigr) \bigr\} \Bigr),$$

*and*

$$\Delta_s = \int \mathbf{I}\bigl\{ x \in W \bigr\}\, \mu_s(dx) - \int \mathbf{I}\bigl\{ x \in W \bigr\}\, \nu_s(dx). \tag{2.45}$$

*If* $\omega(2\sigma_0) < 1/4$ *then*

$$v \le \operatorname{vol}\Bigl( W \cap \bigl\{ 1 - 3\,\omega(2\sigma_0) \le M(x/R) \le 1 + 3\,\omega(2\sigma_0) \bigr\} \Bigr). \tag{2.46}$$

*Proof.* We omit the elementary proof of (2.46).

Using $s_1 = R/\bigl(1+\omega(2\sigma_0)\bigr) \ge R/2$ and the fact that $\omega$ is a nondecreasing function, it is easy to see that

$$s_1\bigl(1+\omega(\sigma_1)\bigr) \le R \le s_2\bigl(1-\omega(\sigma_2)\bigr), \quad \text{with } \sigma_1 = \frac{kr+1}{s_1}, \quad \sigma_2 = \frac{kr+1}{s_2}. \tag{2.47}$$

Let us prove (2.44). Assume first that $\Delta^\circ \ge 0$. Applying (2.20)–(2.22) to $\mu_s$ and $\nu_s$ and using (2.47), we have

$$\begin{aligned} &\mu_{s_2}(C) = p_2 \operatorname{vol}_{\mathbb{Z}} C, \quad \text{for } C \subset R\,\Omega, \\ &\nu_{s_2}(C) \le p_2 \operatorname{vol}\Bigl( C \cap \bigl( s_2(1+\omega(2\sigma_0))\,\Omega \bigr) \Bigr), \quad \text{for any} \quad C. \end{aligned} \tag{2.48}$$

Using (2.48), we have

$$\begin{aligned} &\operatorname{vol}_{\mathbb{Z}}\bigl(W \cap (R\Omega)\bigr) = p_2^{-1}\,\mu_{s_2}\bigl(W \cap (R\Omega)\bigr) \le p_2^{-1} \int \mathbf{I}\bigl\{ x \in W \bigr\}\, \mu_{s_2}(dx), \\ &p_2^{-1} \int \mathbf{I}\bigl\{ x \in W \bigr\}\, \nu_{s_2}(dx) \le \operatorname{vol}\bigl(W \cap (R\Omega)\bigr) + v, \end{aligned} \tag{2.49}$$

proving the lemma in the case $\Delta^\circ \ge 0$.

Assume now that $\Delta^\circ < 0$. Write $p_1 = 1/\operatorname{vol}_{\mathbb{Z}}(s_1\Omega)$. Relations (2.20)–(2.22) together with (2.47) yield

$$\begin{aligned} &\operatorname{vol}\bigl(W \cap (R\Omega)\bigr) \le \operatorname{vol}\Bigl( W \cap \bigl( s_1(1-\omega(2\sigma_0))\Omega \bigr) \Bigr) + v, \\ &\operatorname{vol}\bigl(W \cap (s_1(1-\omega(2\sigma_0))\Omega)\bigr) \le p_1^{-1} \int \mathbf{I}\bigl\{ x \in W \bigr\}\, \nu_{s_1}(dx), \\ &\operatorname{vol}_{\mathbb{Z}}\bigl(W \cap (R\Omega)\bigr) \ge p_1^{-1} \int \mathbf{I}\bigl\{ x \in W \bigr\}\, \mu_{s_1}(dx), \end{aligned}$$

and using $p_2 \le p_1$ we obtain (2.44). □

*Proof of Theorem* 2.6. The result follows from the bound

$$\Delta^* \ll_{d,\varepsilon,m,p,k} q^{d-1}\Big(1 + \tfrac{1}{r^2}\,\tfrac{R^d}{r^d}\Big)\big(\omega(2\sigma_0) + \varepsilon_4 + \varepsilon_5\big) + \frac{q^{d+p}}{\beta-\alpha}\Big(\frac{R^2}{r^2 T} + \frac{R^{2p+2}}{r^{3p}} + \gamma^{1-8/d-\varepsilon}\big(r^2, T\big)\, T^{\varepsilon}\,\frac{R^2}{r^2}\Big) \tag{2.50}$$

choosing the maximal $p < d/2$, $k = 2p+2$ and using the notation (2.38). Recall, that we assume that the constant $c_2$ in (2.38) is sufficiently small. In particular, this assumption guarantees that $\omega(2\sigma_0)$ is as small as will be required in the auxiliary lemmas used below.

Let us prove (2.50). Dividing the bound (2.44) of Lemma 2.7 by $v_1 \overset{\text{def}}{=} \operatorname{vol}\big(W \cap (R\,\Omega)\big)$, we obtain

$$|\Delta^*| \ll \frac{1}{p_2\, v_1} \max_{s=s_1,s_2} |\Delta_s| + \frac{v}{v_1}\,.$$

Hence, in order to prove (2.50) it suffices to show that

$$p_2^{-1} = \operatorname{vol}_{\mathbb{Z}}(s_2\,\Omega) \ll_d R^d, \tag{2.51}$$

$$v_1 \gg_{d,m} (\beta-\alpha)\, q^{-d/2}\, R^{d-2}, \tag{2.52}$$

$$v \ll_{d,m} (\beta-\alpha)\,\big(\omega(2\sigma_0) + \varepsilon_4 + \varepsilon_5\big)\, q^{(d-2)/2}\, R^{d-2}, \tag{2.53}$$

and, for $s = s_1, s_2$,

$$|\Delta_s| \ll |\Delta_{s,1}| + |\Delta_{s,2}| \tag{2.54}$$

with $\Delta_{s,1}$ and $\Delta_{s,2}$ defined below by (2.57) such that

$$|\Delta_{s,1}| \ll_{d,k,\varepsilon} \frac{q^{d/2}}{r^2 T} + \frac{R^{2p}}{r^{3p}}\, q^{p+d/2} + \gamma^{1-8/d-\varepsilon}\big(r^2, T\big)\, T^{\varepsilon}\, \frac{q^{d/2}}{r^2}\,. \tag{2.55}$$

$$|\Delta_{s,2}| \ll_{d,m} (\beta-\alpha)\,\big(\omega(2\sigma_0) + \varepsilon_4 + \varepsilon_5\big)\, q^{(d-2)/2}\, r^{-4}\,(R/r)^{d-2}. \tag{2.56}$$

To prove (2.51) it suffices to notice that $s_2 \le 2R$ since $\Omega \subset B(1)$ and we assume that $\omega(2\sigma_0)$ is sufficiently small.

The estimate (2.52) follows from (1.13) and (8.10) of Lemma 8.2 choosing $\lambda = 1$ and using the estimate $|a_0|/R \le \sqrt{q}\,|a|/R = \varepsilon_1 \le 1/(2m)$ which is guaranteed by the assumption that the constant $c_2$ is sufficiently small.

The bound (2.53) follows from (2.46) and Lemma 8.3 with $\delta = 3\omega(2\sigma_0) \le 1/4$ which is fulfilled since $c_2$ is small. Applying Lemma 3.8 we use the estimate $\varepsilon_{1,0} \le \varepsilon_1$.

Let us prove (2.54). Let $F(I;R)$ and $F_j(I;R)$ denote the functions $F(I)$ and $F_j(I)$ defined by (2.28) with the underlying parameter $R$ which enters into the definitions of the measures $\mu = \mu_R$ and $\nu = \nu_R$. Then (2.54) holds with

$$\Delta_{s,1} = F(I;s) - F_0(I;s) - \Delta_{s,2}, \qquad \Delta_{s,2} = \sum_{j\in 2\mathbb{N},\ j<p} F_j(I;s). \tag{2.57}$$

For the estimation of $\Delta_{s,1}$ we shall apply Theorem 2.2. Our choice of the constant $c_2$ yields $\omega(2\sigma_0) \le 1/2$. Therefore $s = s_1, s_2$ is asymptotically equivalent to $R$. The functions $s \to F(I;s)$ and $s \to F_j(I;s)$ are differences of the corresponding distribution functions. Furthermore, $|a| \ll R$. Hence, Theorem 2.2 yields (2.55).

It remains to prove (2.56). Let us estimate $F_j(I;s)$. Using $s \asymp R$, $r \asymp \overline{r}$ and (7.4), we have

$$\bigl| F_j(I;s) \bigr| \ll_{d,p} v_2 \, r^{-j-d},$$

where

$$v_2 = \mathrm{vol}\Bigl( W \cap \Bigl\{ x \in \mathbb{R}^d : \ \bigl| M(x/s) - 1 \bigr| \le \omega(c_k \, r/R) \Bigr\} \Bigr).$$

By Lemma 8.3 and a suitable choice of $c_2$, the volume $v_2$ allows the same upper bound as the volume $v$ in (2.53) since $s \asymp R$. Hence, summing these bounds over $j$, we get (2.56). □

## 3. Proof of Theorems 2.1 and 2.2

We shall use the following approximate (see, e.g., relation (8.4) in [BG4]) and precise (see, e.g., Chung [Ch]) formulas for the Fourier-Stieltjes inversion. For any $T > 0$ and any distribution function $F$ of a normalized nonnegative measure with the Fourier-Stieltjes transform

$$\widehat{F}(t) = \int \mathrm{e}\bigl\{ tx \bigr\} \, dF(x), \qquad t \in \mathbb{R},$$

we have

$$F(x) = \tfrac{1}{2} + \tfrac{i}{2\pi} \, \mathrm{V.P.} \int_{-T}^{T} \mathrm{e}\bigl\{ -xt \bigr\} \, \widehat{F}(t) \, \tfrac{dt}{t} + \widetilde{\mathcal{R}} \tag{3.1}$$

with remainder term $\widetilde{\mathcal{R}}$ such that

$$|\widetilde{\mathcal{R}}| \le \tfrac{1}{T} \int_{-T}^{T} |\widehat{F}(t)| \, dt.$$

Here $\mathrm{V.P.} \int f(t) \, dt = \lim_{\tau \to 0} \int_{|t| > \tau} f(t) \, dt$ denotes the Principal Value of the integral. Furthermore, for any function $F : \mathbb{R} \to \mathbb{R}$ of bounded variation such that $F(-\infty) = 0$ and $2F(x) = F(x+) + F(x-)$, for all $x \in \mathbb{R}$, we have

$$F(x) = \tfrac{1}{2} F(\infty) + \tfrac{i}{2\pi} \lim_{M \to \infty} \mathrm{V.P.} \int_{|t| \le M} \mathrm{e}\{ -xt \} \, \widehat{F}(t) \, \tfrac{dt}{t} \, . \tag{3.2}$$

The formula is well-known for distribution functions. To functions of bounded variation it extends by linearity arguments.

Theorem 2.1 is implied by Theorem 2.2 and we have to prove Theorem 2.2 only.

The expansion (2.8) yields

$$|\mathcal{R}| = \Big| F - \sum_{j\in 2\,\mathbb{N}_0,\ j<p} F_j \Big|. \tag{3.3}$$

Let us prove Theorem 2.2 assuming that $r \le 1$. Using (3.3), Lemma 7.4 to bound $|F_j|$, and the obvious estimates $|F| \le 1$ and $|F_0| \le 1$, we obtain

$$|\mathcal{R}| \ll_d 1 + \sum_{1\le j<p} \tfrac{R^j}{r^{2j}} \Big(1 + \tfrac{|a|}{r}\Big)^j q^{j+d/2} \ll_d \tfrac{R^p}{r^{2p}} \Big(1 + \tfrac{|a|}{r}\Big)^p q^{p+d/2} \tag{3.4}$$

since $q \ge 1$, $1 \le R/r$, $1/r \ge 1$ and $j < p < d/2$. The estimate (3.4) implies the theorem in the case $r \le 1$. Therefore in the remaining part of the proof we can and shall assume that $r \ge 1$.

Using the representation (3.3), representing $F$ by the approximate Fourier-Stieltjes inversion (3.1) and $F_j$ by the inversion formula (3.2), splitting the intervals of integration, and using the triangle inequality and the obvious estimate

$$\tfrac{1}{T} \int_{1/r\le|t|\le T} \big|\widehat{F}(t)\big|\, dt \le \int_{1/r\le|t|\le T} \big|\widehat{F}(t)\big|\, \tfrac{dt}{|t|},$$

we obtain

$$|\mathcal{R}| \ll I_1 + I_2 + I_3 + \sum_{j\in 2\,\mathbb{N}_0,\ j<p} I_{j+4} \tag{3.5}$$

with

$$I_1 = \int_{|t|\,r\le 1} \Big|\widehat{F}(t) - \sum_{j\in 2\,\mathbb{N}_0,\ j<p} \widehat{F}_j(t)\Big|\, \tfrac{dt}{|t|},$$

$$I_2 = \tfrac{1}{T} \int_{|t|\,r\le 1} \big|\widehat{F}(t)\big|\, dt,$$

$$I_3 = \int_{|t|\,r\ge 1,\ |t|\le T} \big|\widehat{F}(t)\big|\, \tfrac{dt}{|t|},$$

$$I_{j+4} = \int_{|t|\,r\ge 1} \big|\widehat{F}_j(t)\big|\, \tfrac{dt}{|t|}, \qquad j \in 2\,\mathbb{N}_0,\ j < p.$$

The estimate (3.5) shows that in order to prove the theorem it suffices to prove that

$$I_1,\, I_{j+4} \ll_{d,k} \tfrac{R^p}{r^{2p}} \Big(1 + \tfrac{|a|}{r}\Big)^p q^{p+d/2}, \qquad j \in 2\,\mathbb{N}_0,\ j < p, \tag{3.6}$$

$$I_2 \ll_d \tfrac{q^{d/2}}{r^2 T}, \tag{3.7}$$

$$I_3 \ll_{d,\varepsilon} \gamma^{1-8/d-\varepsilon}\gamma\big(r^2,T\big)\,T^{\varepsilon}\,\frac{q^{d/2}}{r^2}\,, \quad \text{for any } \varepsilon>0. \tag{3.8}$$

Let us estimate $I_3$. Changing variables, we have

$$\widehat{F}(t)=\int \mathrm{e}\big\{\,tQ[x-a]\,\big\}\,\mu(dx).$$

Splitting

$$\mu=\mu^{*3}(\cdot\,;r)*\chi \quad \text{with } \chi=\Phi*\mu^{*(k-3)}(\cdot\,;r), \tag{3.9}$$

we obtain

$$I_3\le I_3', \quad I_3' \overset{\text{def}}{=} \sup_{a\in\mathbb{R}^d}\int_{r\,|t|\ge 1,\ |t|\le T} \varphi_a(t;r^2)\,\frac{dt}{|t|}\,, \tag{3.10}$$

where $\varphi_a$ is given by

$$\varphi_a(t;r^2)=\Big|\int \mathrm{e}\big\{\,tQ[x-a]\,\big\}\,\mu^{*3}(dx;r)\Big|, \tag{3.11}$$

(cf. (2.6)). The function $\varphi_a$ satisfies the following inequalities (see Theorem 5.1 in [BG4])

$$\varphi_a(t;r^2)\,\varphi_a(t+\tau;r^2)\ll_d q^{d/2}\,\mathcal{M}^{d/2}(\tau;r^2), \tag{3.12}$$

$$\varphi_a(\tau;r^2)\ll_d q^{d/2}\,\mathcal{M}^{d/2}(\tau;r^2), \quad t,\,\tau\in\mathbb{R} \tag{3.13}$$

involving the function $\mathcal{M}(t;r^2)$ introduced in (1.21). The inequality (3.12) allows us to apply Theorem 5.1 of the present paper. Choosing in this theorem

$$\Lambda=c_d\,q^{d/2}, \quad \varkappa=d/2, \quad s=r^2, \quad \alpha=-1,$$

and using $\ln x\ll_\varepsilon x^\varepsilon$, for $x\ge 1$ and $\varepsilon>0$, we obtain

$$I_3'\ll_{d,\varepsilon}\left(\frac{\gamma\Big(r^2,T\Big)}{q^{d/2}}\right)^{1-8/d-\varepsilon}(1+\ln T)\,\frac{q^{d/2}}{r^2}\,. \tag{3.14}$$

Estimating in (3.14) $1+\ln T\ll_\varepsilon T^\varepsilon$, for $T\ge 1$, and $q\ge 1$, and using (3.10), we derive the desired bound (3.8) for $I_3$.

Let us estimate $I_2$. Similarly to the estimation of $I_3$ we obtain

$$I_2\le \frac{I_2'}{T}\,, \qquad I_2' \overset{\text{def}}{=} \sup_{a\in\mathbb{R}^d}\int_{r\,|t|\le 1}\varphi_a(t;r^2)\,dt. \tag{3.15}$$

The bound (3.13) for $\varphi_a$ and the definition (1.21) of the function $\mathcal{M}(t;r^2)$ together with the inequality $\varphi_a\le 1$ yield

$$\varphi_a(t;r^2)\ll_d q^{d/2}\min\big\{\,1;\ (r^2\,|t|)^{-d/2}\,\big\}, \qquad \text{for } r\,|t|\le 1. \tag{3.16}$$

Hence, we have

$$I_2'\ll_d q^{d/2}\int_{\mathbb{R}}\min\big\{\,1;\ (r^2\,|t|)^{-d/2}\,\big\}\,dt=\frac{q^{d/2}}{r^2}\int_{\mathbb{R}}\min\big\{\,1;\ |t|^{-d/2}\,\big\}\,dt\ll_d \frac{q^{d/2}}{r^2}\,,$$

which together with (3.15) yields the bound (3.7) for $I_2$.

Let us prove (3.6) for $I_1$. Applying the bound (4.10) to $|\mathcal{R}|$ defined by (3.3), bounding by $\varphi_a$ the function $\psi$ in (4.10) and using (3.16), we obtain

$$I_1 \ll_{d,k} \frac{q^{p+d/2}R^p}{r^{2p}} \Big(1+\frac{|a|}{r}\Big)^p I_1'$$

with

$$I_1' = \int\limits_{r\,|t|\le 1} \Big(\big(|t|r^2\big)^{\theta(p)} + \big(|t|r^2\big)^p\Big) \min\big\{1;\ (|t|r^2)^{-d/2}\big\} \frac{dt}{|t|}. \tag{3.17}$$

By change of variables $t r^2 = \tau$ we get

$$I_1' \ll \int\limits_0^\infty (\tau^{\theta(p)} + \tau^p) \min\big\{1;\ \tau^{-d/2}\big\} \frac{d\tau}{\tau} \ll_{d,p} 1, \tag{3.18}$$

since $0 < \theta(p) \le p < d/2$, see (4.7). The bounds (3.17) and (3.18) yield (3.6) for $I_1$.

Let us prove (3.6) for $I_{j+4}$. Using (7.13) we have

$$I_{j+4} \ll_d \frac{R^j}{r^{j+d/2}} \Big(1+\frac{|a|}{r}\Big)^j q^{j+d/2};$$

an application of the inequalities

$$r \ge 1, \quad R/r \ge 1, \quad j < p < d/2, \quad q \ge 1$$

implies (3.6) for $I_{j+4}$, thus concluding the proof of Theorem 2.2.

## 4. An expansion for the Fourier-Stieltjes transform $\widehat{F}$

By change of variables we can write the Fourier-Stieltjes transforms of the distribution functions $F$ and $F_0$ as

$$\widehat{F}(t) = \int \mathrm{e}\big\{ tQ[x-a]\big\}\,\mu(dx), \qquad \widehat{F_0}(t) = \int \mathrm{e}\big\{ tQ[x-a]\big\}\,\nu(dx). \tag{4.1}$$

Similarly, for the Fourier-Stieltjes transforms of the functions of bounded variation $F_j$ (see (2.13)–(2.15)) we have

$$\widehat{F}_j(t) = \sum\nolimits^{**}_{\eta:\ |\eta|_1=j} \widehat{F}_{j\eta}(t), \tag{4.2}$$

where the sum $\sum\nolimits^{**}_{\eta:\ |\eta|_1=j}$ is defined as in (2.13), and

(4.3)
$$\widehat{F}_{j\eta}(t) = \frac{(-1)^m}{\eta!} \int \cdots \int \mathrm{e}\big\{ tQ[x-a]\big\}\, D^{(j)}(x) u_1^{\eta_1} \dots u_m^{\eta_m}\, dx \prod_{l=1}^m \pi^{*(k+1)}(du_l).$$

Introducing the notation

$$g(x) = \exp\big\{ h(x)\big\}, \qquad h(x) = itQ[x-a], \tag{4.4}$$

integrating by parts and using $D(x)\,dx = \nu(dx)$, we have

$$\widehat{F}(t) = \int g(x)\,\mu(dx), \qquad \widehat{F_0}(t) = \int g(x)\,\nu(dx) \tag{4.5}$$

and

$$\widehat{F}_{j\eta}(t) = \frac{(-1)^m}{\eta!} \int\cdots\int g^{(j)}(x)u_1^{\eta_1}\dots u_m^{\eta_m}\,\nu(dx)\prod_{l=1}^m \pi^{*(k+1)}(du_l). \tag{4.6}$$

In the proof of (4.6) we used that the function $D$ has compact support and the derivatives $\partial^\alpha D$ are continuous functions, for $|\alpha|_\infty \le k-2$.

Write

$$\theta(s) = s/2, \text{ for even } s, \qquad \text{and} \qquad \theta(s) = (s+1)/2, \text{ for odd } s, \tag{4.7}$$

and

$$\psi(t) = \sup_{a\in\mathbb{R}^d}\Big|\int \mathrm{e}\big\{tQ[x] + \langle a,x\rangle\big\}\,\mu^{*[k/(p+1)]}(dx;r)\Big|, \tag{4.8}$$

$$\psi_0(t) = \sup_{a\in\mathbb{R}^d}\Big|\int \mathrm{e}\big\{tQ[x] + \langle a,x\rangle\big\}\,\nu^{*[k/(p+1)]}(dx;\overline{r})\Big|. \tag{4.9}$$

LEMMA 4.1. *Let $2 \le p \le k-1$ and $q \ge 1$. Then, for $R \ge r \ge 1$, we have*

$$\Big|\widehat{F}(t) - \widehat{F_0}(t) - \sum_{j\in 2\mathbb{N},\ j<p}\widehat{F_j}(t)\Big| \ll_{d,k} \frac{\psi(t)\,R^p\,q^p}{r^{2p}}\Big(1+\frac{|a|}{r}\Big)^p\Big(\big(|t|r^2\big)^{\theta(p)} + \big(|t|r^2\big)^p\Big) \tag{4.10}$$

*and, for $j \in 2\mathbb{N}_0$, $j < p$,*

$$\big|\widehat{F_j}(t)\big| \ll_{d,k} \frac{\psi_0(t)\,R^j\,q^j}{r^{2j}}\Big(1+\frac{|a|}{r}\Big)^j\Big(\big(|t|r^2\big)^{\theta(j)} + \big(|t|r^2\big)^j\Big). \tag{4.11}$$

*Proof.* Let $g : \mathbb{R}^d \to \mathbb{C}$ denote a sufficiently smooth complex valued function. Assume that $x, u_1, \dots, u_p \in \mathbb{R}^d$. We shall use the following decomposition

$$g(x) = g(x+u_1) + g_1 + \cdots + g_p, \tag{4.12}$$

where

$$g_j = \sum\nolimits^*_{\eta:\ |\eta|_1=j} c(\eta)\,g^{(j)}(x+u_{m+1})u_1^{\eta_1}\dots u_m^{\eta_m}, \quad \text{for } 1 \le j < p, \tag{4.13}$$

$$g_p = \sum\nolimits^*_{\eta:\ |\eta|_1=p} c_p(\eta)\int_0^1 (1-\tau)^{\eta_m}\,g^{(p)}(x+\tau u_m)u_1^{\eta_1}\dots u_m^{\eta_m}\,d\tau, \tag{4.14}$$

and

$$c(\eta) = \frac{(-1)^m}{\eta_1!\dots\eta_m!}, \qquad c_p(\eta) = \frac{(-1)^m}{\eta_1!\dots\eta_{m-1}!\,(\eta_m-1)!}. \tag{4.15}$$

Here the sum $\sum\nolimits^*_{\eta:\ |\eta|_1=j}$ is taken over all representations of $j$ as a sum $j = \eta_1+\cdots+\eta_m$, $m \le j$, of integers $\eta_1, \dots, \eta_m \in \mathbb{N}$. For example, in the case $j = 3$ we have the following four representations: $3=3$, $3=2+1$, $3=1+2$, $3=1+1+1$. In particular, if $p=1$ then

$$g_1 = -\int_0^1 (1-\tau)\, g'(x+\tau u_1)u_1\, d\tau,$$

and if $p=2$ then

$$g_1 = -g'(x+u_2)u_1,$$
$$g_2 = \int_0^1 (1-\tau)\, g''(x+\tau u_2)u_1u_2\, d\tau - \int_0^1 (1-\tau)^2\, g''(x+\tau u_1)u_1^2\, d\tau.$$

In order to prove (4.12) it suffices to iteratively apply Taylor expansions. In the first step we use the Taylor expansion of the function $g$:

(4.16)
$$g(x) = g(x+u_1) - \sum_{\eta_1=1}^{p-1} \frac{1}{\eta_1!}\, g^{(\eta_1)}(x)u_1^{\eta_1} - \frac{1}{(p-1)!}\int_0^1 (1-\tau)^p\, g^{(p)}(x+\tau u_1)u_1^p\, d\tau.$$

In the second step we apply expansions of type (4.16) to the functions $x \to g^{(\eta_1)}(x)u_1^{\eta_1}$, for $1 \le \eta_1 < p$, using $u_2$ instead of $u_1$. After $p$ such steps we arrive at (4.12).

For the derivatives of the composite function $g(x)$ (see (4.4)), we shall use the following formula, which follows from a general formula (see Theorem 2.5 in Averbuh and Smoljanov [AS]). Let $v_1, \dots, v_s \in \mathbb{R}^d$. Let $\mathcal{A} = \{A_1, \dots, A_\alpha\}$ denote a partition of the set $\{v_1, \dots, v_s\} = A_1 \cup \dots \cup A_\alpha$ into $\alpha$ disjoint subsets $A_j$ such that $1 \le \operatorname{card} A_j \le 2$, for all $j$. Notice that $\theta(s) \le \alpha \le s$, where $\theta(s)$ is defined by (4.7). For a subset $A \subset \{v_1, \dots, v_s\}$, say $A = \{v_{i_1}, \dots, v_{i_l}\}$, introduce the derivative $\partial_A h(x) \overset{\text{def}}{=} h^{(l)}(x)v_{i_1}\dots v_{i_l}$. Then

$$g^{(s)}(x)v_1\dots v_s = g(x)\sum_{\mathcal{A}} \partial_{A_1}h(x)\dots\partial_{A_\alpha}h(x), \tag{4.17}$$

where the sum $\sum_{\mathcal{A}}$ extends over all partitions $\mathcal{A}$ with properties specified above. One can easily prove (4.17) using $h'''(x) \equiv 0$ and induction in $s$.

The following identities

$$\nu(\cdot;\overline{r}) = \mu(\cdot;r) * \pi, \qquad \nu = \mu * \pi^{*(k+1)}, \tag{4.18}$$

are obvious (see (2.2), (2.3) and (2.4) for the definitions of measures which appear in (4.18)). For example, for any integrable function $u : \mathbb{R}^d \to \mathbb{R}$, we have

$$\int u(x)\,\nu(dx;\overline{r}) = \sum_{y\in\mathbb{Z}^d\cap B(r)} (2\overline{r})^{-d} \int u(x+y)\,\pi(dx) = \iint u(x+y)\,\pi(dx)\,\mu(dy;r),$$

proving the first identity in (4.18).

Integrating both sides of the identity (4.12) with respect to the measures

$$\mu(dx),\ \ \pi^{*(k+1)}(du_1),\ \ \dots,\ \ \pi^{*(k+1)}(du_p)$$

and using (4.5), (4.6), (4.13), (4.14), (4.18), we obtain

$$\widehat{F}(t) = \widehat{F}_0(t) + f_1(t) + \dots + f_{p-1}(t) + f_p(t), \tag{4.19}$$

with $\widehat{F}$ and $\widehat{F}_0$ defined by (4.5),

$$f_j(t) = \sum_{\eta:\ |\eta|_1=j}^{*} c(\eta) \int \cdots \int g^{(j)}(x) u_1^{\eta_1} \dots u_m^{\eta_m}\, \nu(dx) \prod_{s=1}^{m} \pi^{*(k+1)}(du_s), \tag{4.20}$$

for $1 \le j < p$, and

$$f_p(t) = \sum_{\eta:\ |\eta|_1=p}^{*} c_p(\eta) \int_0^1 (1-\tau)^{\eta_m} f_p(t;\eta)\, d\tau, \tag{4.21}$$

$$f_p(t;\eta) \overset{\text{def}}{=} \int \cdots \int g^{(p)}(x+\tau u_m) u_1^{\eta_1} \dots u_m^{\eta_m}\, \mu(dx) \prod_{l=1}^{m} \pi^{*(k+1)}(du_l).$$

The sums in (4.20) and (4.21) are the same as in (4.13) and (4.14) respectively. However, notice that $f_j = 0$, for odd $j < p$, since the measure $\pi$ is symmetric. Moreover, $f_j = \widehat{F}_j$, for even $j < p$, since all terms in the sum (4.20) vanish unless all $\eta_1, \dots, \eta_m$ are even, due to the symmetry of $\pi$. Hence, (4.19) yields

$$\Big| \widehat{F}(t) - \widehat{F}_0(t) - \sum_{j \in 2\mathbb{N},\ j<p} \widehat{F}_j(t) \Big| = \big| f_p(t) \big|. \tag{4.22}$$

Now we can return to the proof of (4.10). The equality (4.22) shows that it suffices to verify that $\big| f_p(t) \big|$ is bounded from above by the right-hand side of (4.10). Using (4.21), we have to verify that any of the suprema $\sup_{0\le\tau\le1} \big| f_p(t;\eta) \big|$ is bounded by the right-hand side of (4.10), for all allowable $\eta_1, \dots, \eta_s$. Define $v_1, \dots, v_p$ repeating $u_1$ $\eta_1$ times, followed by $u_2$ $\eta_2$ times, etc. For example, $v_j = u_1$, for all $1 \le j \le \eta_1$. Using (4.17) and (4.21) with $m = p$, we obtain

$$\sup_{0\le\tau\le1} \big| f_p(t;\eta) \big| \ll_p \sup_{0\le\tau\le1} \max_{\mathcal{A}} \int \cdots \int \sup_{|z|_\infty \le 1} I_{\mathcal{A}} \prod_{l=1}^{m} \pi^{*(k+1)}(du_l), \tag{4.23}$$

where

$$I_{\mathcal{A}} = \Big| \int g(x+z)\, \partial_{A_1} h(x+z) \dots \partial_{A_\alpha} h(x+z)\, \mu(dx) \Big|, \qquad z \overset{\text{def}}{=} \tau u_m, \tag{4.24}$$

and $\theta(p) \le \alpha \le p$. Fix a partition $\mathcal{A} = \big\{ A_1, \dots, A_\alpha \big\}$ of type used in (4.17) and (4.23). Let $\beta$ denote the number of 1-point sets in this partition. Without loss of generality we can assume that $A_1, \dots, A_\beta$ are 1-point sets, and that $A_{\beta+1}, \dots, A_\alpha$ are 2-point sets. Then

$$\partial_{A_j} h(x+z) = 2it \big\langle x - a + z, Q w_j \big\rangle, \qquad \text{for } 1 \le j \le \beta, \tag{4.25}$$

and

$$\partial_{A_j} h(x+z) = 2it\langle \overline{w}_j, Qw_j\rangle, \qquad \text{for } \beta < j \le \alpha, \tag{4.26}$$

with some $w_j, \overline{w}_j \in \{u_1, \dots, u_m\}$. Notice that $|w_j|_\infty \le 1$ and $|\overline{w}_j|_\infty \le 1$. Furthermore, using (4.25), (4.26) and (4.24), substituting

$$\langle x - a + z, Qw_j\rangle = \langle x, Qw_j\rangle - \langle a - z, Qw_j\rangle,$$

multiplying, applying the triangle inequality and re-enumerating $w_j$ if necessary, we obtain

$$I_{\mathcal{A}} \ll_p |t|^\alpha \max_{0\le\rho\le\beta} I_\rho \tag{4.27}$$

where

$$I_\rho = \Big| \int g(x+z) \prod_{j=1}^{\rho} \langle x, Qw_j\rangle \, B \, \mu(dx)\Big|, \tag{4.28}$$

and

$$B \stackrel{\text{def}}{=} \prod_{j=\rho+1}^{\beta} \langle z - a, Qw_j\rangle \prod_{j=\beta+1}^{\alpha} \langle \overline{w}_j, Qw_j\rangle. \tag{4.29}$$

Using $|Qw| \le q\,|w| \ll_d q$, we have

$$|B| \ll_{d,p} q^{\alpha-\rho}\big(|z|^{\beta-\rho} + |a|^{\beta-\rho}\big) \ll_p q^{\alpha-\rho}\big(1 + |a|^{\beta-\rho}\big) \tag{4.30}$$

since $|z|_\infty = \tau\,|u_m|_\infty \le 1$ (see (4.24)) and $|w_j|_\infty \ll 1$, $|\overline{w}_j|_\infty \ll 1$.

Let us split the measure $\mu = \Phi * \mu^{*k}(\cdot\,; r)$ as follows

$$\mu = \chi_0 * \chi^{*(p+1)}, \quad \chi_0 = \Phi * \mu(\cdot\,; r)^{*(k-(p+1)\,[k/(p+1)])}, \quad \chi = \mu^{*[k/(p+1)]}(\cdot\,; r), \tag{4.31}$$

where $[u]$ denotes the integer part of $u$. Then we have

$$\int U(x)\,\mu(dx) = \int U(x_0 + \dots + x_{p+1})\,\chi_0(dx_0)\,\chi(dx_1)\dots\chi(dx_{p+1}), \tag{4.32}$$

for any integrable function $U$. Applying (4.32) to (4.28), writing $x = x_0 + \dots + x_{p+1}$ and using

$$\prod_{j=1}^{\rho} \langle x, Qw_j\rangle = \sum_{j_1=0}^{p+1} \cdots \sum_{j_\rho=0}^{p+1} \langle x_{j_1}, Qw_1\rangle \dots \langle x_{j_\rho}, Qw_\rho\rangle,$$

we obtain

$$I_\rho \ll_p \max_{0\le j_1,\dots,j_\rho\le p+1} I_{\rho j}, \tag{4.33}$$

with

$$I_{\rho j} = \int \Big| \int \cdots \int g(x+z)\,\langle x_{j_1}, Qw_1\rangle \dots \langle x_{j_\rho}, Qw_\rho\rangle\, B \prod_{s=1}^{p+1} \chi(dx_s)\Big|\, \chi_0(dx_0). \tag{4.34}$$

Given the variables $x_{j_1}, \dots, x_{j_\rho}$, $\rho \le p$, we find among $x_1, \dots, x_{p+1}$ at least one variable, say $x_l$, such that $l \notin \{ j_1, \dots, j_\rho \}$. Without loss of generality we can assume that $l = 1$. Then (4.34) yields

$$I_{\rho j} \ll \int \cdots \int \Big| \big\langle Q x_{j_1}, w_1 \big\rangle \dots \big\langle Q x_{j_\rho}, w_\rho \big\rangle \Big| \, |B| \, J \, \chi_0(dx_0) \prod_{s=2}^{p+1} \chi(dx_s) \tag{4.35}$$

with

$$J = \Big| \int g(x+z) \, \chi(dx_1) \Big|, \qquad x = x_0 + \cdots + x_{p+1}.$$

Recall (see (4.4)) that $g(x) = \mathrm{e}\big\{ tQ[x-a] \big\}$. Therefore

$$J \le \sup_{L \in \mathbb{R}^d} \Big| \int \mathrm{e}\big\{ tQ[x] + \big\langle L, x \big\rangle \big\} \, \mu^{*[k/(p+1)]}(dx; r) \Big| = \psi(t) \tag{4.36}$$

with $\psi$ defined by (4.8). Hence, the bound (4.35) combined with (4.30) and (4.36) yields

$$I_{\rho j} \ll_{d,p} \psi(t) \, q^{\alpha - \rho} \big( 1 + |a|^{\beta - \rho} \big) \, J_{\rho j} \tag{4.37}$$

with

$$J_{\rho j} = \int \cdots \int \Big| \big\langle Q x_{j_1}, w_1 \big\rangle \dots \big\langle Q x_{j_\rho}, w_\rho \big\rangle \Big| \, \chi_0(dx_0) \prod_{s=2}^{p+1} \chi(dx_s). \tag{4.38}$$

We have

$$\Big| \big\langle Q x_{j_1}, w_1 \big\rangle \dots \big\langle Q x_{j_\rho}, w_\rho \big\rangle \Big| \ll_k q^\rho \, R^\rho \tag{4.39}$$

since we assume that $R \ge r$, and since the variables $x_s$ in the integral (4.38) satisfy $|x_s|_\infty \ll_k r + R \ll R$. Combining (4.37)–(4.39), we get

$$I_{\rho j} \ll_k \psi(t) \, R^\rho \, q^\alpha \big( 1 + |a|^{\beta - \rho} \big). \tag{4.40}$$

The estimate (4.40) combined with (4.33) and (4.27) yields

$$I_{\mathcal{A}} \ll \psi(t) \, \big( |t| r^2 \big)^\alpha \, q^\alpha \, r^{-2\alpha} \max_{0 \le \rho \le \beta} R^\rho \, \big( 1 + |a|^{\beta - \rho} \big). \tag{4.41}$$

Using the condition $R \ge r \ge 1$ and $\beta + 2(\alpha - \beta) = 2\alpha - \beta = p$, $\rho \le \beta$, we obtain

(4.42)

$$r^{-2\alpha} \max_{0 \le \rho \le \beta} R^\rho \, \big( 1 + |a|^{\beta - \rho} \big) \ll r^{-2\alpha} \max_{0 \le \rho \le \beta} R^\rho \, r^{\beta - \rho} \Big( 1 + \Big( \tfrac{|a|}{r} \Big)^{\beta - \rho} \Big)$$
$$\ll_p r^{-p} \, \tfrac{R^\beta}{r^\beta} \, \Big( 1 + \tfrac{|a|}{r} \Big)^\beta.$$

The inequalities $\beta \le p$, $q \ge 1$ and $\theta(p) \le \alpha \le p$ combined with (4.40)–(4.42) and (4.21)–(4.23) yield (4.10).

The proof of (4.11) repeats the proof of (4.10) starting from (4.22) since now we have to estimate $|\widehat{F}_j| = |f_j|$ (see (4.20)) instead of $|f_p|$. In this proof

we have to use (4.20) instead of (4.21) and to replace everywhere $p$ by $j$. Furthermore, instead of (4.31) we have to use a similar splitting of $\nu$ of the form $\nu = \chi_0 * \chi^{*(p+1)}$ with

$$\chi_0 = \Psi * \nu(\,\cdot\,;\overline{r})^{*(k-(p+1)\,[k/(p+1)])} \quad \text{and} \quad \chi = \nu^{*[k/(p+1)]}(\,\cdot\,;\overline{r}). \tag{4.43}$$

In particular, the splitting (4.43) yields that the integral corresponding to $J$ in (4.36) now satisfies

$$J \le \sup_{L\in\mathbb{R}^d} \Big| \int \mathrm{e}\big\{ tQ[x] + \langle L, x\rangle \big\}\, \nu^{*[k/(p+1)]}(dx;\overline{r}) \Big| = \psi_0(t),$$

which leads to the factor $\psi_0$ in (4.11). □

## 5. The integration procedure for large $|t|$

Recall that

$$\mathcal{M}(t;s) = \big(|t|s\big)^{-1}\, \mathbf{I}\big\{ |t| \le s^{-1/2} \big\} + |t|\, \mathbf{I}\big\{ |t| > s^{-1/2} \big\}, \tag{5.1}$$

where $s > 0$ will be a positive large parameter.

For a number $T \ge 1$ and a family of functions $\varphi(\cdot) = \varphi(\cdot\,; s) : \mathbb{R} \to \mathbb{R}$ introduce

$$\gamma = \gamma\big(s, T\big) \stackrel{\text{def}}{=} \sup\Big\{ \big|\varphi(t)\big| :\ s^{-1/2} \le t \le T \Big\}. \tag{5.2}$$

The following Theorem 5.1 sharpens Theorem 6.1 of [BG4] in cases where $\gamma < 1$.

Theorem 5.1. *Let $\varphi(t)$, $t \ge 0$ denote a continuous function such that $\varphi(0) = 1$ and $0 \le \varphi \le 1$. Assume that, for some $\varkappa > 4$ and $\Lambda \ge 1$,*

$$\varphi(t)\, \varphi(t+\varepsilon) \le \Lambda \mathcal{M}^{\varkappa}(\varepsilon; s), \quad \textit{for all } t \ge 0 \quad \textit{and} \quad \varepsilon \ge 0. \tag{5.3}$$

*Let $T \ge 1$. Assume that the number $\gamma$ defined by (5.2) satisfies*

(5.4)
$$\gamma > 4^{\varkappa/(\varkappa-4)}\, s^{-\varkappa/4}, \quad \textit{if } \ -1 < \alpha \le 0,$$
$$\gamma \left(1 + \ln \tfrac{1}{\gamma}\right)^{\varkappa/(\varkappa-4)} > 4^{\varkappa/(\varkappa-4)}\, s^{-\varkappa/4}\, (1 + \ln s)^{\varkappa/(\varkappa-4)}, \quad \textit{if } \ \alpha = -1.$$

*Then the integral*

$$J = \int_{s^{-1/2}}^{T} \varphi(t)\, t^{\alpha}\, dt$$

*can be bounded as follows:*

$$J \ll_{\alpha,\varkappa} \left(\frac{\gamma}{\Lambda}\right)^{1-4/\varkappa} T^{\alpha+1}\, \frac{\Lambda}{s}, \qquad \textit{for } \ -1 < \alpha \le 0, \tag{5.5}$$

*and*

$$J \ll_{\varkappa} \left(\frac{\gamma}{\Lambda}\right)^{1-4/\varkappa} \left(1+\ln\frac{\Lambda}{\gamma}\right)(1+\ln T)\,\frac{\Lambda}{s}, \qquad \textit{for } \alpha=-1.$$

*If* (5.4) *is not fulfilled then the following trivial bounds hold*

(5.6)
$$J \ll_{\alpha} \gamma T^{\alpha+1}, \quad \textit{for } \alpha>-1, \qquad J \ll \gamma(1+\ln s)(1+\ln T), \quad \textit{for } \alpha=-1.$$

*Proof.* Evaluating the integral $J$ using the method of Lebesgue integration by partitioning the range of $\varphi$ in intervals $[2^{-l-1}, 2^{-l}]$, we have to estimate the Lebesgue measure of the corresponding sets

$$B_l = \{t : 2^{-l-1} \le \varphi(t) \le 2^{-l}\} \cap [s^{-1/2}, T].$$

Using inequality (5.3), we shall show that two points in $B_l$ are either very "close" or far apart. This means that the set $B_l$ consists of 'small' clusters of size $\mathcal{O}_l(s^{-1})$ separated by "large" gaps of size $\mathcal{O}_l(1)$. These constraints on the structure of $B_l$ suffice to bound the measure of $B_l$ well enough in order to estimate the size of $J$ for $\varkappa > 4$ as claimed in Theorem 5.1. For trigonometric sums $\varphi$ this condition translates to to the assumption that the dimension $d$ satisfies $d > 8$.

Throughout the proof we shall write $\ll$ instead of $\ll_{\alpha,\varkappa}$.

To prove (5.6) it suffices to use $\varphi(t) \le \gamma$ and to notice that

$$\int_{s^{-1/2}}^{T} \frac{dt}{t} \ll (1+\ln s)(1+\ln T), \qquad \int_{s^{-1/2}}^{T} t^{\alpha}dt \ll T^{\alpha+1}, \quad \text{for } \alpha>-1.$$

Let us prove (5.5). The inequality (5.3) implies that (set $t=0$, use $\varphi(0)=1$ and note that $\Lambda \ge 1$)

$$\varphi(t) \le \Lambda \mathcal{M}^{\varkappa}(t;s) \quad \text{and} \quad \varphi(t)\,\varphi(t+\varepsilon) \le \Lambda^2 \mathcal{M}^{\varkappa}(\varepsilon;s). \tag{5.7}$$

Starting the proof of (5.5) with (5.7) we may assume without loss of generality that $\Lambda = 1$, that is, that

$$\varphi(t) \le \mathcal{M}^{\varkappa}(t;s) \quad \text{and} \quad \varphi(t)\,\varphi(t+\varepsilon) \le \mathcal{M}^{\varkappa}(\varepsilon;s). \tag{5.8}$$

Indeed, we may replace $\varphi$ in (5.7) (resp. $\gamma$) by $\varphi/\Lambda$ (resp. by $\gamma/\Lambda$), and we may integrate over $\varphi/\Lambda$ instead of $\varphi$. Notice that now $\varphi(0) \le 1$ and the case $\varphi(0) < 1$ is not excluded.

Thus assuming (5.8) we have to prove that

$$\int_{s^{-1/2}}^{T} \varphi(t)\, t^{\alpha}\, dt \ll \gamma^{1-4/\varkappa}\, \frac{F_\alpha}{s}, \quad \text{for} \quad \varkappa > 4, \tag{5.9}$$

with $F_\alpha = T^{\alpha+1}$, for $-1 < \alpha \le 0$, and $F_{-1} = (1+\ln T)\left(1+\ln\frac{1}{\gamma}\right)$. While proving (5.9) we may assume that $1 \ll s$. Otherwise (5.9) obviously holds since $\varphi \le \gamma \le 1$.

Let $l_\gamma$ denote the smallest integer such that $2^{-l_\gamma} \ge \gamma$. For the integers $l \ge l_\gamma$, introduce the sets

$$B_l = [s^{-1/2}, T] \cap \{ t : \ 2^{-l-1} \le \varphi(t) \le 2^{-l} \},$$
$$D_l = [s^{-1/2}, T] \cap \{ t : \ \varphi(t) \le 2^{-l-1} \}.$$

Since the function $\varphi$ satisfies $0 \le \varphi(t) \le \gamma$, for $s^{-1/2} \le t \le T$, the sets $B_l$ and $D_l$ are closed and $D_m \cup \bigcup_{l=l_\gamma}^{m} B_l = [s^{-1/2}, T]$. Furthermore, (5.8) implies that $\varphi(t) \le t^{\varkappa}$, for $t \ge s^{-1/2}$, whence $B_l \subset [L_l^{-1}, T]$, where $L_l = 2^{(l+1)/\varkappa}$.

Recall that $\varphi(t) \le 2^{-l}$, for $t \in B_l$, and $\varphi(t) \le 2^{-m-1}$, for $t \in D_m$. Therefore the relation $D_m \subset [s^{-1/2}, T]$ yields

$$\int_{s^{-1/2}}^{T} \varphi(t)\, t^{\alpha}\, dt \le \int_{D_m} \varphi(t)\, t^{\alpha}\, dt + \sum_{l=l_\gamma}^{m} \int_{B_l} \varphi(t)\, t^{\alpha}\, dt \ll 2^{-m} G_\alpha + \sum_{l=l_\gamma}^{m} 2^{-l} \int_{B_l} t^{\alpha}\, dt, \tag{5.10}$$

where $G_\alpha = \ln T + \ln s$, for $\alpha = -1$, and $G_\alpha = T^{\alpha+1}$, for $-1 < \alpha \le 0$. We shall choose $m$ such that

$$2^{-m} G_\alpha \le \gamma^{1-4/\varkappa} F_\alpha\, s^{-1}, \qquad \text{for } -1 \le \alpha \le 0. \tag{5.11}$$

More precisely, we choose the minimal $m$ such that

$$m \ge \frac{1}{\ln 2} \ \ln \frac{s\, G_\alpha}{\gamma^{1-4/\varkappa} F_\alpha}, \qquad \text{for } -1 \le \alpha \le 0. \tag{5.12}$$

Using (5.10), (5.11), we see that the estimate (5.9) follows provided that we show that

$$\sum_{l=l_\gamma}^{m} I_l \ll \gamma^{1-4/\varkappa}\, \frac{F_\alpha}{s}, \quad \text{where} \quad I_l = 2^{-l} \int_{B_l} t^{\alpha}\, dt. \tag{5.13}$$

Below we shall prove the inequalities

$$\begin{aligned} I_l &\ll (l + \ln T)\, s^{-1}\, 2^{-l+4l/\varkappa}, && \text{for} \quad \alpha = -1, \\ I_l &\ll T^{\alpha+1}\, s^{-1}\, 2^{-l+4l/\varkappa}, && \text{for} \quad -1 < \alpha \le 0, \end{aligned} \tag{5.14}$$

for $l \le m$. These inequalities imply (5.13). Indeed, in both cases $\alpha = -1$ and $\alpha > -1$ we can apply the bound

$$\sum_{l=l_\gamma}^{\infty} 2^{-l+4l/\varkappa} \ll (2^{-l_\gamma})^{1-4/\varkappa} \ll \gamma^{1-4/\varkappa}$$

since $\varkappa > 4$ ensures the convergence of the series, and, according to the definitions of $\gamma$ and $l_\gamma$, we have $2^{-l_\gamma-1} \le \gamma$. In the case $\alpha = -1$ one needs in addition the following estimates. For $l_\gamma \ll 1$, we have

$$\sum_{l=l_\gamma}^{\infty} l\, 2^{-l+4l/\varkappa} \le \sum_{l=0}^{\infty} l\, 2^{-l+4l/\varkappa} \ll 1 \ll (2^{-l_\gamma})^{1-4/\varkappa} \ll \gamma^{1-4/\varkappa},$$

and, for $l_\gamma \ge c_\varkappa$ with a sufficiently large constant $c_\varkappa$, we obtain

$$\sum_{l=l_\gamma}^{\infty} l\, 2^{-l+4l/\varkappa} \ll \int_{l_\gamma-1}^{\infty} x\, 2^{-x+4x/\varkappa}\, dx \ll l_\gamma (2^{-l_\gamma})^{1-4/\varkappa} \ll \gamma^{1-4/\varkappa} \Big(1 + \ln \frac{1}{\gamma}\Big).$$

It remains to prove the inequalities (5.14). For the estimation of $I_l$ we need a description of the structure of the sets $B_l$ with $l \le m$. Let $t, t' \in B_l$ denote points such that $t' > t$. The inequality (5.8) and the definition of $B_l$ imply

$$4^{-l-1} \le \mathcal{M}^\varkappa(t'-t; s). \tag{5.15}$$

If $t' - t \le s^{-1/2}$ then by (5.15) and the definition of $\mathcal{M}(\varepsilon; s)$ we get

$$t' - t \le \delta, \quad \text{where} \quad \delta = s^{-1} 4^{(l+1)/\varkappa}. \tag{5.16}$$

If $t' - t \ge s^{-1/2}$ then by (5.15) and the definition of $\mathcal{M}(\varepsilon; s)$ we have

$$t' - t \ge \rho, \quad \text{where} \quad \rho = 4^{-(l+1)/\varkappa}. \tag{5.17}$$

For $\varkappa > 4$ and sufficiently large $s \gg 1$ note that

$$\delta < \rho, \quad \text{provided} \quad l \le m. \tag{5.18}$$

Indeed, using the definitions of $\delta$ and $\rho$, we see that the inequality (5.18) follows from the inequality $4 s^{-\varkappa/4} < 2^{-m}$, which is implied by (5.12), the assumption (5.4) and $s \gg 1$.

The estimate (5.18) implies that either $t' - t \le \delta$ or $t' - t \ge \rho$. Therefore it follows from (5.16)–(5.18) that

$$t \in B_l \Longrightarrow B_l \cap (t+\delta, t+\rho) = \emptyset; \tag{5.19}$$

that is, that in the interval $(t+\delta, t+\rho)$ the function $\varphi$ takes values lying outside of the interval $[2^{-l-1}, 2^{-l}]$.

Let us return to the proof of (5.14). If the set $B_l$ is empty then (5.14) is obviously fulfilled. If $B_l$ is nonempty then define $e_1 = \min\big\{t : \ t \in B_l\big\}$. Choosing $t = e_1$ and using (5.19) we see that the interval $(e_1 + \delta, e_1 + \rho)$ does not intersect $B_l$. Similarly, let $e_2$ denote the smallest $t \ge e_1 + \rho$ such that $t \in B_l$. Then the interval $(e_2 + \delta, e_2 + \rho)$ does not intersect $B_l$. Repeating this procedure we construct a sequence $L_l^{-1} \le e_1 < e_2 < \cdots < e_k \le T$ such that

$$B_l \subset \bigcup_{j=1}^{k} [e_j, e_j + \delta] \quad \text{and} \quad e_{j+1} \ge e_j + \rho. \tag{5.20}$$

The sequence $e_1 < \cdots < e_k$ cannot be infinite. Indeed, due to (5.20) we have

$$T \geq e_k \geq e_1 + (k-1)\rho \geq L_l^{-1} + (k-1)\,\rho \geq k\,\rho,$$

and therefore $k \leq T/\rho$.

Using (5.20) we can finally prove (5.14). We start with the case $\alpha = -1$. Using $\ln(1+x) \leq x$, for $x \geq 0$, we have

$$\begin{aligned} I_l &\leq 2^{-l} \sum_{j=1}^{k} \int_{e_j}^{e_j+\delta} \frac{dt}{t} = 2^{-l} \sum_{j=1}^{k} \ln\Big\{1 + \frac{\delta}{e_j}\Big\} \\ &\leq 2^{-l} \sum_{j=1}^{k} \frac{\delta}{e_j} \ll (l + \ln T)\, 2^{-l+4l/\varkappa}\, s^{-1} \end{aligned}$$

since $e_1 \geq L_l^{-1} \geq \rho$, $k \leq T/\rho$, and

$$\sum_{j=1}^{k} \frac{1}{e_j} \leq \sum_{j=1}^{k} \frac{1}{e_1 + (j-1)\rho} \leq \frac{1}{\rho} \sum_{j=1}^{k} \frac{1}{j} \ll \frac{\ln T + \ln \rho^{-1}}{\rho} \ll (l + \ln T)\, 4^{l/\varkappa}.$$

Finally, let us prove (5.14) for $-1 < \alpha \leq 0$. We have

$$\begin{aligned} I_l &\leq \frac{1}{2^l} \sum_{j=1}^{k} \int_{e_j}^{e_j+\delta} t^{\alpha} dt \\ &\leq \frac{1}{2^l} \sum_{j=1}^{k} \delta\, e_j^{\alpha} \leq \frac{\delta\,\rho^{\alpha}}{2^l} \sum_{j=1}^{k} j^{\alpha} \ll \frac{\delta\, T^{\alpha+1}}{2^l\,\rho} \ll s^{-1} 2^{-l+4l/\varkappa} T^{\alpha+1}. \end{aligned}$$

□

## 6. Trigonometric sums of irrational quadratic forms

In this section we introduce a criterion for the rationality of a quadratic form in terms of the behavior of an associated sequence of trigonometric sums. Recall that a quadratic form $Q[x] = \langle Qx, x \rangle$, $x \in \mathbb{R}^d$, with a nonzero symmetric matrix $Q = (q_{ij})$, $1 \leq i, j \leq d$, is rational if there exists an $M \in \mathbb{R}$, $M \neq 0$, such that the matrix $MQ$ has integer entries; otherwise it is irrational.

For $a \in \mathbb{R}^d$, $d \geq 1$, consider the polynomial

$$P(x) = Q[x] + \langle a, x \rangle, \qquad x \in \mathbb{R}^d. \tag{6.1}$$

Throughout this section we shall denote by $\mu$ the uniform lattice measure in the box $B(r) = \big\{ |x|_\infty \leq r \big\}$ (cf. with the notation $\mu(\cdot\,; r)$ used in other sections). In other words, the measure $\mu$ is nonnegative, normalized $\mu(\mathbb{R}^d) = \mu\big(\mathbb{Z}^d \cap B(r)\big) = 1$ and assigns equal weights $\mu_x \overset{\text{def}}{=} \mu(\{x\}) = \big(2\,[r] + 1\big)^{-d}$ to points $x \in \mathbb{Z}^d \cap B(r)$.

For $k \in \mathbb{N}$ and $t \in \mathbb{R}$, introduce the trigonometric sum

$$f(t) = f(t;r) = \sum_{x\in\mathbb{Z}^d} \mathrm{e}\bigl\{ t P(x) \bigr\}\, \mu_x^{*(2k+1)}, \quad \mu_x^{*(2k+1)} \overset{\mathrm{def}}{=} \mu^{*(2k+1)}(\{x\}),$$

where $\mu^{*k}$ denotes the $k$-fold convolution of $\mu$. Let $\widetilde{\mu}$ denote the symmetrization of $\mu$, that is, $\widetilde{\mu}(C) = \int \mu(C+x)\,\mu(dx)$, for $C \in \mathcal{B}^d$. Writing $\mu^{*(2k+1)} = \mu * \mu^{*k} * \mu^{*k}$ and using the symmetrization inequality (see Lemma 6.5 below), we obtain $0 \le \bigl| f(t)\bigr|^2 \le \varphi(t)$ with

$$\varphi(t) = \varphi(t;r) \overset{\mathrm{def}}{=} \sum_{x\in\mathbb{Z}^d} \sum_{y\in\mathbb{Z}^d} \mathrm{e}\bigl\{ 2t\langle Qx, y\rangle \bigr\}\, \widetilde{\mu}_x\, \widetilde{\mu}_y^{*k}. \tag{6.2}$$

For a family of functions $g = g(\cdot\,;r) : \mathbb{R} \to \mathbb{C}$ parameterized by $r \ge 0$, consider the following condition:

$$\text{for any } \ 0 < \delta_0 \le \delta < \infty, \qquad \lim_{r\to\infty} \sup_{\delta_0 \le |t| \le \delta} \bigl| g(t;r)\bigr| = 0. \tag{6.3}$$

We shall apply condition (6.3) to the trigonometric sums $f$ and $\varphi$. See (6.4) and (6.6) for equivalent formulations of (6.3). The following characterization result holds without assumptions on the eigenvalues of $Q$.

THEOREM 6.1. *Let $k \ge 1$. The quadratic form $Q[x]$ is irrational if and only if $\varphi$ satisfies condition* (6.3). *If $Q$ is irrational then $f$ satisfies* (6.3). *If $a$ and $Q$ are rational* (*that is*, *there exists $M \ne 0$ such that $Ma$ and $MQ$ have rational coordinates*, *resp. entries*) *then $f$ does not satisfy* (6.3).

The proof of Theorem 6.1 will be given later. It is based on an application of the theory of successive minima (see Cassels [Ca], Davenport [Dav]) and techniques in [BG1].

For a given family of functions, the condition (6.3) allows the following equivalent reformulation: there exist functions $\delta_0(r) \downarrow 0$ and $\delta(r) \uparrow \infty$ such that

$$\lim_{r\to\infty} \sup_{\delta_0(r) \le |t| \le \delta(r)} \bigl| g(t;r)\bigr| = 0. \tag{6.4}$$

Applying a double large sieve type bound (see the estimate (5.22) in [BG4]) we obtain $\varphi(t) \ll_{d,k} q^{2d} \mathcal{M}^d(t; r^2)$ with $\mathcal{M}$ defined by (1.21). Hence, assuming $\delta_0(r) \ge r^{-1}$, we have $\varphi(t) \ll_{d,k,Q} \delta_0^d(r)$, for $r^{-1} \le t \le \delta_0(r)$, and

$$\lim_{r\to\infty} \sup_{a\in\mathbb{R}^d} \sup_{r^{-1} \le t \le \delta_0(r)} \varphi(t) = 0 \tag{6.5}$$

provided that the eigenvalues of $Q$ are nonzero. Combining (6.4) and (6.5) we see that the irrationality of $Q$ is equivalent to the following condition: there exist $\delta(r) \uparrow \infty$ such that

$$\lim_{r\to\infty} \sup_{a\in\mathbb{R}^d} \sup_{r^{-1}\le t\le\delta(r)} \varphi(t) = 0 \tag{6.6}$$

provided that the eigenvalues of $Q$ are nonzero. Due to $|f|^2 \le \varphi$, relations (6.5) and (6.6) hold for $f$ as well.

*Remark* 6.2. An inspection of the proof shows that the condition (6.3) for the trigonometric sums $f$ and $\varphi$ holds uniformly over compact sets, say $\mathcal{Q}$, of irrational matrices; that is,

$$\lim_{r\to\infty} \sup_{Q\in\mathcal{Q}} \sup_{\delta_0\le|t|\le\delta} \bigl|\varphi(t;r)\bigr| = 0.$$

Let us recall some facts of the theory of successive minima in the geometry of numbers (see [Dav]). Let $F : \mathbb{R}^d \to [0,\infty)$ be a norm, that is, $F(\alpha x) = |\alpha| F(x)$, for $\alpha \in \mathbb{R}$, and $F(x+y) \le F(x) + F(y)$. The successive minima $M_1 \le \cdots \le M_d$ of $F$ with respect to the lattice $\mathbb{Z}^d$ are defined as follows: $M_1 = \inf\bigl\{F(x) :\ x \ne 0,\ x \in \mathbb{Z}^d\bigr\}$, and $M_k$ is defined as the lower bound of $\lambda > 0$ such that the set $\bigl\{x \in \mathbb{Z}^d :\ F(x) < \lambda\bigr\}$ contains $k$ linearly independent vectors. It is easy to see that there exist linearly independent vectors $a_1, \dots, a_d \in \mathbb{Z}^d$ such that $F(a_j) = M_j$.

LEMMA 6.3 ([Dav, Lemma 3]). *Let $L_j(x) = \sum_{k=1}^d \lambda_{jk} x_k$, $1 \le j \le d$, be linear forms in $\mathbb{R}^d$ such that $\lambda_{jk} = \lambda_{kj}$. Assume that $P \ge 1$ and let $\|\theta\|$ denote the distance of the number $\theta$ to the nearest integer. Then the number of $x = (x_1, \dots, x_d) \in \mathbb{Z}^d$ such that*

$$\|L_j(x)\| < P^{-1}, \qquad |x|_\infty < P, \qquad \text{for all } 1 \le j \le d,$$

*is bounded from above by $\frac{c_d}{M_1 \cdots M_d}$, where $M_1 \le \cdots \le M_d$ are the first $d$ of the $2d$ successive minima $M_1 \le \cdots \le M_{2d}$ of the norm $F : \mathbb{R}^{2d} \to [0,\infty)$ defined for $y = (x,m) \in \mathbb{R}^{2d}$ with $x, m \in \mathbb{R}^d$ as*

$$F(y) \overset{\text{def}}{=} \max\Bigl\{P\bigl|L_1(x) - m_1\bigr|, \dots, P\bigl|L_d(x) - m_d\bigr|,\ P^{-1}|x|_\infty\Bigr\}. \tag{6.7}$$

*Furthermore, $M_1 \ge P^{-1}$.*

We shall use some well-known known bounds for trigonometric sums which go back to Weyl [We]. For our purposes we formulate a variation of such bounds as the following Lemma 6.4. Its proof is provided for the sake of completeness.

LEMMA 6.4. *For* $z=(z_1,\dots,z_d)\in\mathbb{R}^d$,

$$g(z)\overset{\text{def}}{=}\sum_{y\in\mathbb{Z}^d}\mathrm{e}\big\{\langle z,y\rangle\big\}\,\widetilde{\mu}_y^{*k}\ll_{d,k}\sum_{m\in\mathbb{Z}^d}h(m),$$

*where*

$$h(m)=\prod_{j=1}^{d}h_0(z_j-2\pi m_j),\qquad h_0(s)=\big(1+r^2s^2\big)^{-k},\qquad s\in\mathbb{Z}. \tag{6.8}$$

*Proof.* We have $g(z)=v^k(z)$ with

$$v(z)=\sum_{y\in\mathbb{Z}^d}\mathrm{e}\big\{\langle z,y\rangle\big\}\,\widetilde{\mu}_y=\Big|\sum_{y\in\mathbb{Z}^d}\mathrm{e}\big\{\langle z,y\rangle\big\}\,\mu_y\Big|^2=\prod_{j=1}^{d}u(z_j), \tag{6.9}$$

where $u(x)=D_{[r]}^2(x)/\big(2\overline{r}\big)^2$ and $2\overline{r}=2[r]+1$ with the Dirichlet kernel

$$D_{[r]}(x)=\sum_{n=-[r]}^{[r]}\exp\big\{ixn\big\}=\frac{\sin(x\overline{r})}{\sin(x/2)},\qquad\text{for }\ x\in\mathbb{R}.$$

Using $|\sin y|\le\min\big\{1;\,|y|\big\}$ and $x^2\ll\sin^2(x/2)$, with $y=x\overline{r}$ and $|x|\le\pi$, we obtain

$$u(x)\ll\min\Big\{\frac{1}{x^2\,\overline{r}^2}\,;\,1\Big\}\ll\frac{1}{1+r^2\,x^2},\qquad\text{for}\quad|x|\le\pi. \tag{6.10}$$

The function $u(x)$ is even and $2\pi$-periodic. Hence (6.10) yields

$$u^k(x)\ll_k\sum_{m\in\mathbb{Z}}\frac{\mathbf{I}\Big\{|x-2\pi m|\le\pi\Big\}}{\Big(1+r^2\,(x-2\pi m)^2\Big)^k}\ll\sum_{m\in\mathbb{Z}}\frac{1}{\Big(1+r^2\,(x-2\pi m)^2\Big)^k}, \tag{6.11}$$

where $\mathbf{I}\{A\}$ denotes the indicator function of the event $A$. Combining $g(z)=v^k(z)$ and (6.9), (6.11) we conclude the proof. □

*Proof of Theorem* 6.1. Let us show that for rational $Q$ the trigonometric sum $\varphi$ does not satisfy (6.3). Assume for simplicity that $Q$ has integer entries. Let $\mu\times\nu$ denote the product measure of measures $\mu$ and $\nu$. Clearly, $\big\langle Qx,y\big\rangle\in\mathbb{Z}$, for $x,y\in\mathbb{Z}^d$, and the function

$$\varphi(t)=\sum_{m\in\mathbb{Z}}p_m\,e\big\{2tm\big\},\qquad p_m\overset{\text{def}}{=}\widetilde{\mu}\times\widetilde{\mu}^{*k}\big(\big\{(x,y)\in\mathbb{Z}^{2d}:\ \big\langle Qx,y\big\rangle=m\big\}\big),$$

with weights $p_m$ such that $\sum_{m\in\mathbb{Z}}p_m=1$, is a $\pi$-periodic function and $\varphi(\pi)=1$. Similar arguments show that for rational $Q$ and $a$ the trigonometric sum $f$ does not satisfy (6.3).

To complete the proof of the theorem, it remains to prove that if $f$ or $\varphi$ does not satisfy condition (6.3) then $Q$ is rational. The inequality $|f|^2\le\varphi$ shows that we have to consider the case of $\varphi$ only. Since (6.3) is not fulfilled, we have

$$\text{for some } 0<\delta_0\le\delta<\infty, \qquad \limsup_{r\to\infty}\sup_{\delta_0\le|t|\le\delta}\varphi(\pi t;r)>0.$$

Hence, there exist $\varepsilon>0$ and sequences $t_n\to t_0$, $\delta_0\le|t_n|\le\delta$, and $r_n\uparrow\infty$ such that

$$\varphi(\pi t_n;r_n)\ge\varepsilon>0, \qquad \text{for all } n\in\mathbb{N}. \tag{6.12}$$

We shall prove that (6.12) implies the rationality of $Q$.

Henceforth we shall write $t$ and $r$ instead of $t_n$ and $r_n$. Without loss of generality we shall assume throughout the proof that $r=r_n$ is sufficiently large, that is, that $r\ge c$ with some sufficiently large constant $c=c(d,k,Q,\delta,\varepsilon)$. Furthermore, we shall write $\ll$ instead of $\ll_{d,k,Q,\delta,\varepsilon}$.

The measure $\widetilde{\mu}$ is concentrated in the cube $B(2r)$ such that

$$\sum_{x\in\mathbb{Z}^d}G(x)\,\widetilde{\mu}_x=\sum_{x\in B(2r)\cap\mathbb{Z}^d}G(x)\widetilde{\mu}_x, \qquad \sum_{x\in B(2r)\cap\mathbb{Z}^d}\widetilde{\mu}_x=1, \quad 0\le\widetilde{\mu}_x\ll r^{-d}, \tag{6.13}$$

for any summable function $G$. Using the representation (6.2) of $\varphi$, Lemma 6.4, the relations (6.12) and (6.13) for the function $h(m)$ defined by (6.8) with $z=2\pi tQx$, we have

$$\varepsilon\ll r^{-d}\sum_{x\in B(2r)\cap\mathbb{Z}^d}\;\sum_{m\in\mathbb{Z}^d}h(m). \tag{6.14}$$

Obviously we may replace $h(m)$ in (6.14) by

$$h(m)\overset{\text{def}}{=}\prod_{j=1}^{d}h_0(z_j-m_j), \qquad h_0(s)=\big(1+r^2s^2\big)^{-k}, \qquad z=z(x)=tQx, \tag{6.15}$$

since this yields in the inequality (6.14) an extra factor depending on $d$ and $k$ only.

For vectors $s=(s_1,\dots,s_d)\in\mathbb{Z}^d$ such that $s_j\ge0$ and $\rho=(\rho_1,\dots,\rho_d)\in\mathbb{Z}^d$ such that $\rho_j=\pm1$, for all $j$, introduce the class $H_{s,\rho}$ of vectors $(x,m)\in\mathbb{Z}^{2d}$ such that

$$x\in B(2r), \quad \operatorname{sign}\big(z_j(x)-m_j\big)=\rho_j, \quad \tfrac{s_j}{4r}\le|z_j(x)-m_j|<\tfrac{(s_j+1)}{4r},$$

for all $1\le j\le d$. Here $z_j(x)$ denote the $j^{\text{th}}$ coordinate of $z(x)=tQx$. Moreover, consider the class $H$ of vectors $(x,m)\in\mathbb{Z}^{2d}$ such that

$$x\in B(4r), \quad \|z_j(x)\|<\tfrac{1}{4r}, \quad \text{for all } 1\le j\le d, \tag{6.16}$$

where $\|\theta\|$ denotes the distance from the number $\theta$ to the nearest integer. Then

$$\bigcup_{s,\rho}H_{s,\rho}=\big\{(x,m)\in\mathbb{Z}^d:\ x\in B(2r)\big\}.$$

Furthermore, it is easy to see that $(x,m),(\overline{x},\overline{m}) \in H_{s,\rho}$ implies that $(x-\overline{x}, m-\overline{m}) \in H$. Hence, card $H_{s,\rho} \le$ card $H$. Using (6.14) and the definitions introduced, we obtain

$$\begin{aligned}(6.17)\qquad \varepsilon &\ll r^{-d} \sum_{s,\rho} \text{card } H_{s,\rho} \ \max\bigl\{ h(m) : \ (x,m) \in H_{s,\rho} \bigr\} \\ &\ll r^{-d} \text{card } H \sum_{s} \prod_{j=1}^{d} \bigl(1+s_j^2\bigr)^{-k} \ll r^{-d} \text{card } H,\end{aligned}$$

for sufficiently large $r$.

Lemma 6.3 shows that the cardinality of the class $H$ defined by (6.16) is bounded from above by $\frac{c_d}{M_1 \cdots M_d}$, where $M_j$ are the successive minima of the norm $F$ defined by (6.7) with $P=4r$ and

$$(6.18)\qquad L_j(x) = z_j(x) = (tQx)_j = \sum_{i=1}^{d} t q_{ji} x_i.$$

Combining (6.17) with card $H \ll (M_1 \cdots M_d)^{-1}$, we obtain $r^d M_1 \cdots M_d \ll 1$. Hence, using $P=4r$ and $P^{-1} \le M_1 \le \cdots \le M_d$, we have $M_j \asymp r^{-1}$, for all $1 \le j \le d$.

Let $(a_s, b_s) \in \mathbb{Z}^{2d}$ denote linearly independent vectors such that $F\bigl((a_s,b_s)\bigr) = M_s$, for all $1 \le s \le 2d$. The equalities $F\bigl((a_s,b_s)\bigr) = M_s$ together with (6.7) and relations $M_s \asymp r^{-1} \asymp P^{-1}$ imply that

$$(6.19)\qquad |a_s|_\infty \ll 1, \quad |b_s|_\infty \ll 1, \quad 1 \le s \le d.$$

Combining the relations $F\bigl((a_s,b_s)\bigr) = M_s$ with $P=4r$ and (6.7), (6.18), we get (recall that $t=t_n$ and $r=r_n$; therefore $a_s = a_s^{(n)}$ and $b_s = b_s^{(n)}$ also depend on $n$)

$$(6.20)\qquad \Bigl| \sum_{k=1}^{d} t_n q_{jk} a_{sk}^{(n)} - b_{sj}^{(n)} \Bigr| \ll \frac{1}{r_n^2}, \qquad 1 \le j \le d, \quad 1 \le s \le d.$$

The inequalities (6.19) guarantee that the sequence of systems $\bigl\{ a_1^{(n)}, \dots, a_d^{(n)} \bigr\}$ in $\mathbb{Z}^d$ repeats infinitely often. Choosing an appropriate subsequence, we may assume that $a_s^{(n)} = a_s$ are independent of $n$. Similarly we may assume that $b_s^{(n)} = b_s$ are independent of $n$. Passing to the limit as $n \to \infty$ in (6.20) along the subsequence and using $t_n \to t_0$, $t_0 \ne 0$, we see that the numbers $t_0 q_{j1}, \dots, t_0 q_{jd}$ satisfy

$$(6.21)\qquad \sum_{i=1}^{d} a_{si} t_0 q_{ji} = b_{sj}, \qquad 1 \le s \le d,$$

for all $1 \le j \le d$. Below we shall prove that the vectors $a_1, \dots, a_d$ are linearly independent. Therefore the system (6.21) has the unique solution $t_0 q_{j1}, \dots, t_0 q_{jd}$, which obviously has to be *rational* by Cramer's rule, for all $1 \le j \le d$.

To conclude the proof we have to show that the vectors $a_1, \dots, a_d \in \mathbb{Z}^d$ are linearly independent in $\mathbb{Z}^d$ (or, equivalently, in $\mathbb{R}^d$). If $a_1, \dots, a_d$ are linearly dependent then there exist $v_s \in \mathbb{R}$ not all equal to zero such that $|v_s| < 1$ and $\sum_{s=1}^{d} v_s a_s = 0$. Let us prove that there exist integers $u_1, \dots, u_d$ not all equal to zero such that $|u_s| \ll 1$ and $x \overset{\text{def}}{=} \sum_{s=1}^{d} u_s a_s = 0$. By the multivariate Dirichlet approximation (see, for example, [Ca, §V.10]), for any $N \in \mathbb{N}$ there exist $u_j \in \mathbb{Z}$ and an integer $0 < q \le N$ such that

$$\left| v_s - \frac{u_s}{q} \right| < \frac{1}{q N^{1/d}}, \qquad \text{for all } 1 \le s \le d. \tag{6.22}$$

The inequalities $|v_s| < 1$ and (6.22) yield $|u_s| \le 2N$. Since the vectors $a_s$ have integer coordinates and $|a_s|_\infty \ll 1$, the equation $\sum_{s=1}^{d} v_s a_s = 0$ together with (6.22) implies $x = \sum_{s=1}^{d} u_s a_s = 0$, for sufficiently large $N \ll 1$. Hence $|u_s| \ll 1$.

For the vector $(x, m)$ with $m \overset{\text{def}}{=} \sum_{s=1}^{d} u_s b_s$ we have

$$F\big((x,m)\big) \le \sum_{s=1}^{d} |u_s| F\big((a_s, b_s)\big) \ll \frac{1}{r} \tag{6.23}$$

since $|u_s| \ll 1$ and $F\big((a_s, b_s)\big) \asymp 1/r$. Using (6.7), we see that

$$F\big((x,m)\big) \ge 4r \Big| \sum_{s=1}^{d} u_s b_{sj} \Big|, \qquad \text{for all } 1 \le j \le d, \tag{6.24}$$

where $b_s = (b_{s1}, \dots, b_{sd})$. Combining (6.23) and (6.24), and using that $r$ is sufficiently large and that $\sum_{s=1}^{d} u_s b_{sj}$ are integers, we conclude that $\sum_{s=1}^{d} u_s b_{sj} = 0$, for all $1 \le j \le d$. In other words, $m = \sum_{s=1}^{d} u_s b_s = 0$, which together with the assumption $x = \sum_{s=1}^{d} u_s a_s = 0$ means that the vectors $(a_s, b_s) \in \mathbb{Z}^{2d}$, $1 \le s \le d$, are linearly dependent, a contradiction. □

*A symmetrization inequality.* The following symmetrization inequality is a generalization of a well-known classical inequality due to Weyl [We]. For a proof, see [BG4, Lemma 7.1]. Recall that the symmetrization $\widetilde{\mu}$ of a measure $\mu$ is defined by $\widetilde{\mu}(C) = \int \mu(C + x)\, \mu(dx)$, for $C \in \mathcal{B}^d$.

LEMMA 6.5. *Let $Q : \mathbb{R}^d \to \mathbb{R}^d$ be a linear symmetric operator, $L \in \mathbb{R}^d$ and $C \in \mathbb{R}$. Let $\mu_1, \mu_2, \mu_3, \nu$ denote arbitrary probability measures on $\mathbb{R}^d$. Define a real valued polynomial of second order by*

$$P(x) = \langle Qx, x\rangle + \langle L, x\rangle + C, \qquad \textit{for} \ \ x \in \mathbb{R}^d.$$

*Then the integral*

$$J = \Big| \int \mathrm{e}\big\{ t\, P(x) \big\}\, \mu_1 * \mu_2 * \mu_3 * \nu(dx) \Big|^2$$

*satisfies $2J \le J_1 + J_2$, where*

$$J_1 = \iint \mathrm{e}\big\{ 2t \,\langle Qx, y\rangle \big\}\, \widetilde{\mu}_1(dx)\, \widetilde{\mu}_2(dy),$$

$$J_2 = \iint \mathrm{e}\big\{ 2t \,\langle Qx, z\rangle \big\}\, \widetilde{\mu}_1(dx)\, \widetilde{\mu}_3(dz).$$

*In particular, if $\mu_2 = \mu_3$ then $J \le J_1$.*

## 7. Properties of the distribution functions $F_j$ and the signed measures $\nu_j$

We start by establishing some properties of $F_j$ and $\nu_j$ (see (2.13)–(2.15) for definitions) using the Fourier transforms of the measures $\nu_j$. Recall that we denote

$$W = \big\{ x \in \mathbb{R}^d : \ Q[x - a] \in I \big\}, \quad I = (\alpha, \beta], \quad F_j(I) = F_j(\beta) - F_j(\alpha).$$

Let var $F_j$ denote the variation of $F_j$. Write

$$V = \Big\{ x \in \mathbb{R}^d : \ \big| M(x/R) - 1 \big| \le \omega(\sigma_0) \Big\}, \quad \sigma_0 = (kr + 1)/R, \quad \overline{r} = [r] + 1/2.$$

LEMMA 7.1. *The density $D : \mathbb{R}^d \to \mathbb{R}$ defined by* (2.12) *has continuous bounded partial derivatives $\partial^\alpha D$, for $|\alpha|_\infty \le k - 2$.*

*Assume that $j \le k - 2$. Then we have*

$$\big| D_j(x) \big| \ll_{k,d} \overline{r}^{\,-j-d}\, \mathbf{I}\big\{ |x|_\infty \le \overline{R} + k\overline{r} \big\}, \tag{7.1}$$

*and*

$$\big| F_j(s) \big| \le \mathrm{var}\ F_j \ll_{k,d} \overline{r}^{\,-j}\, (\overline{R}/\overline{r})^d, \quad \textit{for} \ \ R \ge r. \tag{7.2}$$

*If the measure $\Phi$ is defined by* (2.17) *then we have in addition*

$$\big| D_j(x) \big| \ll_{k,d} \overline{r}^{\,-j-d}\, \mathbf{I}\big\{ x \in V \big\}, \tag{7.3}$$

*and*

$$\big| F_j(I) \big| \ll_{k,d} \overline{r}^{\,-j-d}\ \mathrm{vol}(W \cap V). \tag{7.4}$$

*Proof.* The Fourier transforms of complex valued functions $f : \mathbb{R}^d \to \mathbb{C}$ are denoted by

$$\widetilde{f}(y) = \int \mathrm{e}\big\{ \langle y, x \rangle \big\} f(x) \, dx.$$

Let $u(x) = \mathbf{I}\big\{ |x|_\infty \le 1/2 \big\}$, $x \in \mathbb{R}^d$, denote density of the measure $\pi$. The measure $\nu(\cdot\,; \overline{r})$ has density $u_{\overline{r}}(x) \overset{\text{def}}{=} (2\overline{r})^{-d} u\big(x/(2\overline{r})\big)$. Write $U = \frac{d\Psi}{dx}$ for density of the measure $\Psi$ defined by (2.11). Then $\nu$ has density $D = U * u_{\overline{r}}^{*k}$ (see (2.12)), where $*$ denotes the convolution of functions. Therefore $\widetilde{D} = \widetilde{U}\ \widetilde{u}_{\overline{r}}^{k}$ and

$$\big|\partial^\alpha D(x)\big| \ll_d \int |y^\alpha| \, \big|\widetilde{D}(y)\big| \, dy \le \int |y^\alpha| \, \big|\widetilde{u}_{\overline{r}}(y)\big|^k \, dy, \tag{7.5}$$

since $|\widetilde{U}| \le 1$. Hence, using

$$\widetilde{u}(y) = \prod_{j=1}^{d} \frac{\sin(y_j/2)}{y_j/2}, \quad \widetilde{u}_{\overline{r}}(y) = \widetilde{u}(2\overline{r}y),$$

we obtain

$$\big|\partial^\alpha D(x)\big| \ll_{k,d} \overline{r}^{\,-|\alpha|_1 - d}, \quad \text{for } |\alpha|_\infty \le k-2, \tag{7.6}$$

which proves the lemma's assertion about $\partial^\alpha D$.

Let us prove (7.1). Using the definition (2.13)–(2.14) of $D_j$, well-known properties of the Fourier transforms, and the equality $j = \eta_1 + \cdots + \eta_m$, estimating in (2.14) $|u_s| \ll_d 1$ and applying (7.5), (7.6), we obtain

$$\begin{aligned} \big|D_j(x)\big| &\ll_j \max_{\eta:\, |\eta|_1 = j} \int \cdots \int \big| D^{(j)}(x) u_1^{\eta_1} \ldots u_m^{\eta_m} \big| \prod_{l=1}^{m} \pi^{*(k+1)}(du_l) \\ &\ll_{j,d} \int \cdots \int \big| \widetilde{D}(y) \langle y, u_1 \rangle^{\eta_1} \cdots \langle y, u_m \rangle^{\eta_m} \big| \, dy \prod_{l=1}^{m} \pi^{*(k+1)}(du_l) \\ &\ll_{k,d} \int |y|^j \, \big|\widetilde{D}(y)\big| \, dy \ll_{j,d} \overline{r}^{\,-j-d}. \end{aligned}$$

Now (7.1) follows since $D(x) = 0$, for $|x|_\infty > \overline{R} + k\overline{r}$.

Let us prove (7.2). Using (2.15) and (7.1), we obtain

$$\big|F_j(s)\big| \le \int \big|D_j(x)\big| \, dx \ll_{k,d} \overline{r}^{\,-j-d} \, (\overline{R} + k\overline{r})^d \ll_{k,d} \overline{r}^{\,-j} \, (\overline{R}/\overline{r})^d.$$

Let us prove (7.3). By (2.21) and (2.22) the density $D$ is constant outside the set $V$. Thus, outside $V$ the derivatives of $D$ vanish and $D_j(x) = 0$, for $x \in \mathbb{R}^d \setminus V$. Thus (7.1) yields (7.3).

Similarly, (7.4) follows from (7.1), (7.3) and the inequality

$$\big|F(I)\big| \le \int \mathbf{I}\Big\{ x \in W \cap V \Big\} \, \big|D_j(x)\big| \, dx. \qquad \square$$

For further investigation of $F_j$ we shall use the following double large sieve type bound (Lemma 7.2 below) which is a useful corollary of the large sieve of Bombieri and Iwaniec [BI]. It follows from Corollary 5.3 in [BG4] replacing in that corollary $q^2$ by $q$.

LEMMA 7.2. *Assume that functions* $g, h : \mathbb{R}^d \to \mathbb{C}$ *satisfy* $\bigl|g(x)\bigr| \le 1$ *and* $\bigl|h(x)\bigr| \le 1$. *Let* $\mu$ *and* $\nu$ *denote arbitrary probability measures on* $\mathbb{R}^d$ *such that*

$$\mu\Bigl(\bigl\{x \in \mathbb{R}^d : \ |x|_\infty \le T\bigr\}\Bigr) = 1 \quad \text{and} \quad \nu\Bigl(\bigl\{x \in \mathbb{R}^d : \ |x|_\infty \le S\bigr\}\Bigr) = 1,$$

*for some* $T > 0$ *and* $S > 0$. *Write*

$$J = \Bigl|\int \Bigl(\int g(x)\, h(y)\, \mathrm{e}\bigl\{t\bigl\langle Qx, y\bigr\rangle\bigr\}\ \mu(dx)\Bigr)\, \nu(dy)\Bigr|^2, \qquad t \in \mathbb{R}.$$

*Then there exists a positive constant* $c_d$, *depending only on the dimension* $d$, *such that*

$$J \ll_d q^d \bigl(1 + ST|t|\bigr)^d \sup_{x\in\mathbb{R}^d} \mu\Bigl(x + \tfrac{c_d B}{|t|\, S}\Bigr) \sup_{x\in\mathbb{R}^d} \nu\Bigl(x + \tfrac{c_d B}{|t|\, T}\Bigr),$$

*where* $B = \bigl\{x \in \mathbb{R}^d : \ |x|_\infty \le 1\bigr\}$.

COROLLARY 7.3. *For an integer* $m \ge 2$, *the function*

$$\psi_0(t) = \sup_{a\in\mathbb{R}^d} \Bigl|\int \mathrm{e}\bigl\{tQ[x] + \bigl\langle a, x\bigr\rangle\bigr\}\, \nu^{*m}(dx; \overline{r})\Bigr|,$$

*satisfies*

$$\psi_0(t) \ll_d q^{d/2} \min\bigl\{1;\ (r^2\,|t|)^{-d/2}\bigr\}. \tag{7.7}$$

*Proof.* We shall derive the result using Lemma 7.2.

If $1 \ge r^2\,|t|$ then (7.7) is obviously fulfilled since $\psi_0(t) \le 1$. Thus we can assume that $r^2\,|t| > 1$.

Using the obvious identity of the type

$$\int f(x)\, \mu_1 * \mu_2 * \mu_3(dx) = \int f(x + y + z)\, \mu_1(dx)\, \mu_2(dy)\, \mu_3(dz),$$

the fact that the form $Q[x]$ is quadratic and Hölder's inequality, we obtain

$$\psi_0^2(t) \le \sup_{a\in\mathbb{R}^d} \int J\, \nu^{*(m-2)}(dz; \overline{r}) \tag{7.8}$$

with

$$J = \Bigl|\int \Bigl(\int g(x)\, h(y)\, \mathrm{e}\bigl\{2t\bigl\langle Qx, y\bigr\rangle\bigr\}\ \nu(dx; \overline{r})\Bigr)\, \nu(dy; \overline{r})\Bigr|^2,$$

where the functions $|g| \le 1$ and $|h| \le 1$ may depend on $z$ and $a$. Choose $S = T = \overline{r}$. The measure $\nu(\cdot\,; \overline{r})$ has nonnegative density (with respect to the Lebesgue measure on $\mathbb{R}^d$) bounded from above by $(2\overline{r})^{-d}$. Therefore $\nu(A;\overline{r}) \le (2\overline{r})^{-d} \operatorname{vol} A$, for any measurable set $A \subset \mathbb{R}^d$. Thus Lemma 7.2 and the assumption $r^2\,|t| > 1$ yield

$$J \ll_d q^d \big(1+\overline{r}^2\,|t|\big)^d\,\overline{r}^{-2d}\,\big(|t|\overline{r}\big)^{-2d} \ll_d \frac{q^d}{r^{2d}\,|t|^d},$$

whence, using (7.8), we derive (7.7). □

LEMMA 7.4. *The distribution function $s \to F_0(s)$ has a bounded continuous derivative for $d \ge 3$. The functions $s \to F_j(s)$ are functions of bounded variation for $2 \le j < d/2$, and*

$$\sup_s \big|F_j(s)\big| \ll_{j,d} \frac{R^j}{r^{2j}}\,\Big(1+\frac{|a|}{r}\Big)^j\,q^{j+d/2},$$

*Moreover, each of the functions $s \to F_j(s)$, $j \in 2\mathbb{N}_0$, has $[(d-1)/2]-j$ bounded continuous derivatives.*

*Proof.* Changing the variables $t = \tau r^{-2}$, we have

$$\int_{\mathbb{R}} \big(|t|r^2\big)^\alpha \min\Big\{1;\ \big(|t|r^2\big)^{-\beta}\Big\}\,dt \ll_{\alpha,\beta} \frac{1}{r^2}, \quad \text{for } -1<\alpha<\beta-1, \tag{7.9}$$

$$\int_{\mathbb{R}} \big(|t|r^2\big)^\alpha \min\Big\{1;\ \big(|t|r^2\big)^{-\beta}\Big\}\,\frac{dt}{|t|} \ll_{\alpha,\beta} 1, \quad \text{for } 0<\alpha<\beta. \tag{7.10}$$

Using $\theta(j) \le j$, the bound (4.11) for $\big|\widehat{F}_j(t)\big|$, the estimate of Corollary 7.3 for $\psi_0(t)$ and (7.9), (7.10), choosing $\beta = d/2$ and $\alpha = \theta(j)$, $\alpha = j$ respectively, we obtain

$$\int_{\mathbb{R}} \big|\widehat{F}_j(t)\big|\,dt \ll_{j,d} \frac{R^j}{r^{2j+2}}\,\Big(1+\frac{|a|}{r}\Big)^j\,q^{j+d/2}, \tag{7.11}$$

for $j \in 2\mathbb{N}_0$ such that $j < -1 + d/2$. Similarly to (7.11), we obtain

$$\int_{\mathbb{R}} \big|\widehat{F}_j(t)\big|\,\frac{dt}{|t|} \ll_{j,d} \frac{R^j}{r^{2j}}\,\Big(1+\frac{|a|}{r}\Big)^j\,q^{j+d/2}, \tag{7.12}$$

for $j \in 2\mathbb{N}$ such that $j < d/2$. Finally we have

$$\int_{r\,|t|\ge 1} \big|\widehat{F}_j(t)\big|\,\frac{dt}{|t|} \ll_{j,d} \frac{R^j}{r^{j+d/2}}\,\Big(1+\frac{|a|}{r}\Big)^j\,q^{j+d/2}, \tag{7.13}$$

for $r \ge 1$ and $j \in 2\mathbb{N}_0$ such that $j < d/2$, and

$$\int_{\mathbb{R}} |t|^\delta\,\big|\widehat{F}_j(t)\big|\,dt < \infty, \qquad \text{for } 0 \le \delta < -j-1+d/2. \tag{7.14}$$

To conclude the proof it suffices to apply the Fourier and Fourier-Stieltjes inversion formulas and to use the bounds (7.11)–(7.14). □

## 8. On the Lebesgue volumes of certain bodies related to indefinite quadratic forms

We shall describe the asymptotic behavior of the volume of the set

$$A = \{ x \in \mathbb{R}^d : \ M(x) \in R\, I_0, \ \ Q[x-a] \in I \}, \qquad R \to \infty, \tag{8.1}$$

where $M$ is the Minkowski functional of the set $\Omega$ defined by (1.13) and $I_0$ and $I = [\alpha, \beta]$ are finite intervals. Throughout the section we assume that $Q[x]$ is an indefinite quadratic form in $\mathbb{R}^d$, $d \ge 3$. Re-enumerating the eigenvalues of $Q$, choosing an appropriate orthonormal basis of $\mathbb{R}^d$ and denoting the coordinates of $x \in \mathbb{R}^d$ in this basis as $\overline{x}_1, \dots, \overline{x}_d$, we may assume that

$$Q[x] = q_1 \overline{x}_1^2 + \dots + q_d \overline{x}_d^2$$

with $q_1, \dots, q_n > 0$, $q_{n+1}, \dots, q_d < 0$ and some $1 \le n \le d/2$. Indeed, if $n > d/2$, we may replace in (8.1) the matrix $Q$ and the interval $I$ by $-Q$ and $-I$ respectively. Let $S = S^{n-1} \times S^{d-n-1}$ denote a direct product of the unit $(n-1)$- and $(d-n-1)$-dimensional spheres

$$r^2 \overset{\text{def}}{=} \overline{x}_1^2 + \dots \overline{x}_n^2 = 1 \quad \text{ and } \quad \rho^2 \overset{\text{def}}{=} \overline{x}_{n+1}^2 + \dots \overline{x}_d^2 = 1$$

with area elements $d\eta_1$ and $d\eta_2$ respectively, and $d\eta \overset{\text{def}}{=} d\eta_1\, d\eta_2$. Write

$$M_0(x) = M\big( \overline{x}_1/\sqrt{|q_1|}, \dots, \overline{x}_d/\sqrt{|q_d|} \big).$$

LEMMA 8.1. *The volume of the set $A$ defined by* (8.1) *satisfies*

$$\lim_{R\to\infty} R^{-d+2} \operatorname{vol} A = |\det Q|^{-1/2}\, \frac{\beta - \alpha}{2} \int_0^\infty u^{d-3} \Big( \int_S \mathbf{I}\{ M_0(u\eta_1 + u\eta_2) \in I_0 \}\, d\eta \Big)\, du. \tag{8.2}$$

*Proof.* We shall use the following representation

$$\operatorname{vol} A = 2^{-1} R^{d-2}\, |\det Q\,|^{-1/2}\, J, \tag{8.3}$$

where

$$J = R^2 \int_0^\infty \Big( \int_{-\infty}^{u^2} \varphi(r,\rho)\, \mathbf{I}\{ \alpha \le R^2 v \le \beta \}\, r^{n-1} \rho^{d-n-2}\, dv \Big)\, du. \tag{8.4}$$

$$\varphi(r,\rho) = \int_S \mathbf{I}\{ M_0(x + a_0/R) \in I_0 \}\, d\eta, \qquad x = r\,\eta_1 + \rho\,\eta_2, \tag{8.5}$$

and $a_0 = \big(\sqrt{|q_1|}\, a_1, \dots, \sqrt{|q_d|}\, a_d\big)$,

(8.6)
$$r = r(u,v) = u, \;\; \rho = \rho(u,v) = \sqrt{u^2 - v}, \qquad 0 \le u < \infty, \;\; -\infty < v \le u^2.$$

For the proof of (8.3) it suffices to write

$$\operatorname{vol} A = \int \mathbf{I}\big\{\, M(x+a) \in R\, I_0 \,\big\}\, \mathbf{I}\big\{\, \alpha \le Q[x] \le \beta \,\big\}\, dx,$$

to change the variables $\overline{x}_j = R y_j / \sqrt{|q_j|}$, $1 \le j \le d$, to use the polar coordinates $r$ and $\rho$ and to perform the change of variables (8.6).

Let us prove (8.2). Assuming that $R \ge |a_0|$, it is easy to see that

$$\mathbf{I}\big\{\, M_0(x + a_0/R) \in I_0 \,\big\} = 0 \quad \text{unless } \; 0 \le r \le c, \;\; 0 \le \rho \le c,$$

with some sufficiently large $c = c(M_0, I_0)$. Consequently, we have $\varphi(r,\rho) = 0$ unless $0 \le u \le c$, $|v| \le c^2$. Therefore, for sufficiently large $R$, we can write

$$J = \int_0^c f(u)\, du \tag{8.7}$$

with

$$f(u) = r^{n-1} R^2 \int_{\alpha/R^2}^{\beta/R^2} \varphi(r,\rho)\, \mathbf{I}\big\{ -\infty \le v \le u^2 \,\big\}\, \rho^{d-n-2}\, dv.$$

The function $\varphi(r,\rho)$ is a continuous function of the variables $u, v$ and $a_0/R$. Therefore, for any $u > 0$, we obtain

$$\lim_{R\to\infty} f(u) = u^{d-3}\, (\beta - \alpha) \int_S \mathbf{I}\big\{ M_0(u\eta_1 + u\eta_2) \in I_0 \,\big\}\, d\eta.$$

The estimate $f(u) \ll_d c^{d-2}\, \big(|\alpha| + |\beta|\big)$ allows us to apply the dominated convergence theorem, and (8.7) yields

$$\lim_{R\to\infty} J = (\beta - \alpha) \int_0^\infty u^{d-3} \int_S \mathbf{I}\big\{ M_0(u\eta_1 + u\eta_2) \in I_0 \,\big\}\, d\eta\, du,$$

which combined with (8.3) concludes the proof of (8.2). □

The relation (2.18) yields

$$(dq)^{-1/2}\, |x| \le M_0(x) \le m\, |x|. \tag{8.8}$$

LEMMA 8.2. *Let* $I_0 = [0, \lambda]$ *and* $\tau = \lambda + |a_0|/R$, $\sigma = \lambda/m - |a_0|/R$. *Then the volume of the set* $A$ *defined by* (8.1) *satisfies*

$$\operatorname{vol} A \ll_d (\beta - \alpha)\, q^{(d-2)/2}\, \tau^{d-2}\, R^{d-2}. \tag{8.9}$$

*If* $\sigma > 0$ *and* $|\alpha| + |\beta| \le \sigma^2 R^2/5$ *then*

$$\operatorname{vol} A \gg_d (\beta - \alpha)\, q^{-d/2}\, \sigma^{d-2}\, R^{d-2}. \tag{8.10}$$

*Proof.* We shall use the representation (8.3)–(8.6). Estimating $1 \le |\det Q| \le q^d$, we see that it suffices to prove that

$$J \ll_d (\beta-\alpha)\, q^{(d-2)/2}\, \tau^{d-2}, \tag{8.11}$$

$$J \gg_d (\beta-\alpha)\, \sigma^{d-2}. \tag{8.12}$$

Let us prove (8.11). Using (8.8), we have

$$M_0(x + a_0/R) \ge (dq)^{-1/2} \big(|x| - |a_0|/R\big)$$

and

$$\varphi(r,\rho) \le \int_S \mathbf{I}\big\{ |x|^2 \le d\,q\tau^2 \big\}\, d\eta \ll_d \mathbf{I}\big\{ u^2 + \rho^2 \le d\,q\,\tau^2 \big\}$$

since $r = u$ and $|x|^2 = r^2 + \rho^2$. Hence, we get

$$J \ll_d R^2\, q^{(d-3)/2}\, \tau^{d-3} \int_0^\infty \Big( \int_{\alpha/R^2}^{\beta/R^2} \mathbf{I}\big\{ u^2 \le d\,q\,\tau^2 \big\}\, dv \Big)\, du \ll_d (\beta-\alpha)\, q^{(d-2)/2}\, \tau^{d-2},$$

proving (8.11).

Let us prove (8.12). Using (8.8), we have $M_0(x+a_0/R) \le m\,\big(|x|+|a_0|/R\big)$ and

$$\varphi(r,\rho) \ge \int_S \mathbf{I}\big\{ |x|^2 \le \sigma^2 \big\}\, d\eta \gg_d \mathbf{I}\big\{ 2\,u^2 + |v| \le \sigma^2 \big\}$$

since $r = u$ and $\rho^2 = u^2 - v \le u^2 + |v|$. Hence, using the condition $|\alpha| + |\beta| \le \sigma^2\, R^2/5$, we obtain

$$\begin{aligned} J &\gg_d R^2 \int_0^\infty \mathbf{I}\Big\{ \frac{\sigma^2}{4} \le u^2 \le \frac{\sigma^2}{3} \Big\} \\ &\quad\cdot \Big( \int_{-\infty}^{u^2} \mathbf{I}\Big\{ |v| \le \frac{\sigma^2}{5} \Big\} \mathbf{I}\big\{ \alpha \le v\,R^2 \le \beta \big\} r^{n-1}\, \rho^{d-n-2}\, dv \Big) du \\ &\gg_d \sigma^{d-3} R^2 \int \mathbf{I}\Big\{ \frac{\sigma^2}{4} \le u^2 \le \frac{\sigma^2}{3} \Big\} \\ &\quad\cdot \Big( \int \mathbf{I}\big\{ \alpha \le v\,R^2 \le \beta \big\}\, dv \Big)\, du \gg_d (\beta-\alpha)\, \sigma^{d-2}, \end{aligned}$$

proving (8.12). □

The definition (2.19) of the modulus of continuity $\omega$ of $M$ yields

$$\big| M_0(x+y) - M_0(x) \big| \le |x|\, \omega\big( d\,\sqrt{q}\, |y|/|x| \big). \tag{8.13}$$

Write

$$\varepsilon_{1,0} = \frac{|a_0|}{R}, \quad \varepsilon_2 = \frac{|\alpha|}{R^2} + \frac{|\beta|}{R^2}, \quad \varepsilon_4 = \omega(2d\,\varepsilon_{1,0}\sqrt{q}), \quad \varepsilon_5 = \omega(30d\,\varepsilon_2\sqrt{q}). \tag{8.14}$$

Lemma 8.3. *Let $I_0 = [1-\delta, 1+\delta]$, $0 \le \delta \le 1/4$. Assume that*

$$\varepsilon_{1,0} \le c, \qquad \varepsilon_2 \le c, \qquad \varepsilon_4 \le c\, q^{-1/2}, \qquad \varepsilon_5 \le c\, q^{-1/2} \tag{8.15}$$

*with some sufficiently small positive constant $c = c(d,m)$ depending on $d$ and $m$ only. Then the volume of the set $A$ defined by* (8.1) *satisfies*

$$\operatorname{vol} A \ll_{d,m} (\beta - \alpha)\,(\delta + \varepsilon_4 + \varepsilon_5)\, R^{d-2}\, q^{(d-2)/2}. \tag{8.16}$$

*Proof.* We shall write $\ll$ instead of $\ll_{d,m}$. Using the representation (8.3)–(8.6) and the inequality $|\det Q\,|^{-1/2} \le 1$, we reduce the proof of (8.16) to the verification of

$$J \ll (\beta - \alpha)\,(\delta + \varepsilon_4 + \varepsilon_5)\, q^{(d-2)/2}. \tag{8.17}$$

Using $|x|^2 = r^2 + \rho^2$, applying the inequalities (8.8) and the triangle inequality, we see that

$$\varphi(r,\rho) = 0 \quad \text{unless } \ 1 - \delta - \varepsilon_{1,0} \ll |x| \ll \sqrt{q}\,(1-\delta) + \varepsilon_{1,0}$$

since now $I_0 = [1-\delta, 1+\delta]$ with $0 \le \delta \le 1/4$. Hence the assumptions (8.15) and $q \ge 1$ imply that

$$\varphi(r,\rho) = 0 \quad \text{unless } \ 1/2 \le |x| \ll \sqrt{q}. \tag{8.18}$$

The variable $v$ in (8.4) satisfies $|v| \le \varepsilon_2 \le c \le 1/8$. Thus, (8.18) shows that we can assume that $r = u$ in (8.4) satisfies $u \ge 1/4$. Using (8.18) we can estimate $r, \rho \le |x| \ll \sqrt{q}$, and we have

$$J \ll q^{(d-3)/2} R^2 \int_{1/4}^{\infty} \Big( \int_{\alpha/R^2}^{\beta/R^2} \varphi(r,\rho)\, dv \Big)\, du. \tag{8.19}$$

Using (8.13) and (8.18), we obtain

$$\Big| M_0\Big(x + \tfrac{a_0}{R}\Big) - M_0(x) \Big| \ll \sqrt{q}\, \omega\Big( d\sqrt{q}\, \tfrac{|a_0|}{R\,|x|} \Big) \le \sqrt{q}\, \varepsilon_4. \tag{8.20}$$

Let us write $x = r\eta_1 + r\eta_2 + (\rho - r)\eta_2$. The condition $|x| \ge 1/2$ yields $r + \rho \ge 1/2$. Therefore, using $u = r \ge 1/4$ and repeating the arguments used for the proof of (8.20), we obtain

$$\tfrac{|\rho - r|}{r} = \tfrac{|v|}{r\,(r+\rho)} \le 8|v| \le 8\varepsilon_2, \quad \big| M_0(x) - u\, M_0(\eta_1 + \eta_2) \big| \ll \sqrt{q}\, \varepsilon_5. \tag{8.21}$$

Using (8.20) and (8.21), we can replace the indicator function in (8.5) by

$$\mathbf{I}\big\{ 1 - \delta - c_0\sqrt{q}\,\varepsilon_4 - c_0\sqrt{q}\,\varepsilon_5 \le u\, M_0(\eta_1 + \eta_2) \le 1 + \delta + c_0\sqrt{q}\,\varepsilon_4 + c_0\sqrt{q}\,\varepsilon_5 \big\}, \tag{8.22}$$

where the constant $c_0 = c_0(d,m)$ depends on $d$ and $m$ only. Now using (8.19), integrating the indicator function (8.22) in the variable $u$, estimating $M_0(\eta_1 + \eta_2) \gg 1$ and integrating in $v$, we obtain (8.17). □

Universität Bielefeld, Bielefeld, Germany
*E-mail addresses*: bentkus@mathematik.uni-bielefeld.de
goetze@mathematik.uni-bielefeld.de